\documentclass[reqno, 11pt]{amsart}
\usepackage{amsmath}
\usepackage{graphicx}
\usepackage{amsfonts}
\usepackage{amssymb}
\usepackage{color}
\usepackage{mathrsfs}
\usepackage{epsfig}
\usepackage{url}
\usepackage{mathrsfs,a4wide}
\usepackage{calligra}

\theoremstyle{plain}

\theoremstyle{definition}
\theoremstyle{remark}
\numberwithin{equation}{section}

\usepackage{amsmath}
\usepackage{amssymb}
\usepackage{graphicx}
\usepackage{amsthm}
\usepackage{thmtools}
\usepackage{enumitem}
\usepackage{hyperref, cleveref}

\usepackage{mathrsfs}

\usepackage{tikz}

\theoremstyle{plain} \declaretheorem[numberwithin = section, name = Theorem,
 refname = {Theorem}, Refname = {Theorem}]{thm}
\theoremstyle{plain} \declaretheorem[numberlike = thm, name = Proposition,
 refname = {Proposition}, Refname = {Proposition}]{prop}
\theoremstyle{plain} \declaretheorem[numberlike = thm, name = Lemma,
refname = {Lemma}, Refname = {Lemma}]{lem}
\theoremstyle{plain} \declaretheorem[numberlike = thm, name = Definition,
 refname = {Definition}, Refname = {Definition}]{df}
\theoremstyle{definition} 
\theoremstyle{definition} \declaretheorem[numberlike = thm, name = Remark,
refname = {Remark}, Refname = {Remark}]{rem}

\DeclareMathOperator {\real}{\mathbb{R}}
\DeclareMathOperator {\BV}{BV}
\DeclareMathOperator{\weak*}{\begin{picture}(10,4)
							\put(0,-2){$\rightharpoonup$}
							\put(3,3){$\ast$}
							\end{picture}}

\DeclareMathOperator {\HN-1}{{\mathcal H}^{N-1}}
\DeclareMathOperator{\gammalimsup}{\Gamma\text{-}\limsup}
\DeclareMathOperator{\gammaliminf}{\Gamma\text{-}\liminf}
\DeclareMathOperator{\gammalim}{\Gamma\text{-}\lim}

\newcommand{\R}{\real}
\newcommand{\e}{\varepsilon}

\title[Surfactants in the two gradient theory of phase transitions]{Surfactants in the two gradient theory of phase transitions}

\author[M. Cicalese]{M. Cicalese}
\address[Marco Cicalese]{Zentrum Mathematik - M7, Technische Universit\"at M\"unchen, Boltzmannstrasse 3, 85748 Garching, Germany}
\email[M. Cicalese]{cicalese@ma.tum.de}
\author[T. Heilmann]{T. Heilmann}
\address[Tim Heilmann]{Zentrum Mathematik - M7, Technische Universit\"at M\"unchen, Boltzmannstrasse 3, 85748 Garching, Germany}
\email[T. Heilmann]{heilmant@ma.tum.de}

\begin{document}

\pagenumbering{arabic}
\pagestyle{headings}

\begin{abstract}

We investigate the influence of surfactants on stabilizing the formation of interfaces in solid-solid phase transitions. The analysis focuses on singularly perturbed van der Waals–Cahn–Hillard-type energies for gradient vector fields, supplemented with a term that accounts for the interaction between the surfactant and the solid. Assuming the potential term to have only two rank-$1$ connected wells, we prove that the effective energy for the formation of an interface decreases when the surfactant segregate to the interface. 
\vskip5pt
\noindent
\end{abstract}

\maketitle
\tableofcontents

\section*{Introduction}

Surfactants (surface active agents) play a pivotal role in influencing phase transitions.  In essence, the primary mechanism driving these effects is the adsorption of surface active molecules onto phase interfaces. This adsorption alters the surface tension, by decreasing the energy penalty associated with the different chemical environments of the different phases. Consequently surfactants exert a profound influence on the stability and morphology of the physical system. The capacity of surfactants to modulate phase transitions has found practical applications in various fields both in fluid-fluid and in solid-solid phase transitions. In the case of solid-solid transitions some examples are provided by crystal growth, metallurgy, and ceramics processing (see \cite{surfactantbook} and the references therein).\\

In this paper we introduce a phase transition model in presence of surfactant working within the framework of the gradient theory of phase transition. More specifically we modify the easiest phase-field model for solid-solid transition introduced in \cite{cfl} in order to account for the interaction between the surfactant and the solid. The model we introduce draws inspiration from the one proposed by Perkins, Sekerka, Warren and Langer for fluid-fluid transitions and analyzed in \cite{fms} (see also \cite{ab, bbm, cfhp} for some extension to more general models of fluid-fluid or multiphase-fluid-fluid phase transitions in presence of surfactants). To fix the ideas, in what follows we first present the latter model. Such a model, motivated by the investigation on foam stability, is a modification of the classical van der Waals-Cahn-Hillard energy functional. More specifically, an integral term accounting for the fluid-surfactant interaction is added to the classical Cahn-Hillard functional, as explained in detail below.

In a given open and bounded set $\Omega\subset\real^N$ (the region occupied by the fluid and the surfactant), one considers a scalar function $u:\Omega\to \real$ and a non-negative function $\rho:\Omega\to [0,+\infty)$ representing the order parameter of the fluid and the density of the surfactant, respectively. As $\varepsilon\to 0$ one is interested in the asymptotic behaviour of the singularly perturbed sequence of energy functionals ${\mathcal E}_\varepsilon:W^{1,2}(\Omega)\times{\mathcal M}(\Omega)\to [0,+\infty)$ defined as 
\begin{equation} \label{intro1}
	{\mathcal E}_\varepsilon(u,\mu) := \int \limits_\Omega \frac{1}{\varepsilon}W(u)
 	+ \varepsilon |\nabla u|^2 + \varepsilon (\rho - |\nabla u|)^2 \,dx,
\end{equation}
where $\mu = \rho \mathcal{L}^N$ denotes the surfactant measure and $W:\real\to[0,+\infty)$ is a double-well potential with wells $\{u:\, W(u)=0\}=\{0,1\}$. The first two terms in the energy define the usual Cahn-Hillard energy functional, namely $$CH_\e(u)=\int \limits_\Omega \frac{1}{\varepsilon}W(u)+ \varepsilon |\nabla u|^2\,dx,$$ which models the energy cost of a phase separation phenomenon in a two-component immiscible fluid. In few words, within this theory, also known as the gradient theory of phase transitions, the phase separation phenomenon corresponds to the transition from the value $0$ to the value $1$ of the order parameter $u$ which represents the local concentration of one of the components of the fluid. The variational limit in the sense of $\Gamma$-convergence (see \cite{b,dm}) of the Cahn-Hillard functional as $\e\to 0$ has been considered by Modica and Mortola in \cite{mod,modmor} (see also \cite{baldo, ft} for analogous results about the asymptotic behaviour of the Cahn-Hillard functional in the case of vector valued order parameters). In \cite{mod} it is proved the pre-compactness in $BV(\Omega;\{0,1\})$ of sequences of phase-fields $u_\e$ with uniformly bounded energy and it is computed the $\Gamma$-limit of $CH_\e$ as $\e\to 0$ with respect to the $L^1$ convergence. Roughly speaking, the limit $u$ of a converging subsequence of $u_\e$ will take only the values $0$ and $1$, partitioning $\Omega$ in the two sets $\{u=0\}$ and $\{u=1\}$ (the two immiscible phases of the fluid) whose common boundary (the phase interface) will correspond to the jump set $S_u$ of the function $u$. Being $u$ of bounded variation, the latter set will have finite ${\mathcal H}^{N-1}$-measure. Up to a multiplicative constant depending on the shape of $W$ the effective asymptotic energy of the system, captured by the $\Gamma$-limit of $CH_\e$, will be proportional to such a perimeter measure. Hence, if we fix the measure of the set $\{u_\e=0\}$ to be strictly smaller than the measure of $\Omega$, both phases will be non empty and the minimal Cahn-Hillard energy as $\e\to 0$ will correspond to the partition of $\Omega$ in the two sets having the least perimeter of the common boundary. Such an energy will be achieved along a sequence $u_\e$ of phase fields with $\nabla u_\e$ concentrating on $S_u$. In this perspective, one can understand the role of the additional third term in ${\mathcal E}_\e(u,\mu)$ which is responsible of the interaction between the surfactant and the fluid. The presence of the additional term will modify the minimizers of $CH_\e$ described above enhancing the phase separation phenomenon to happen in the regions where the surfactant is present. In fact, the last integral term in \eqref{intro1} is minimized if $\rho_\e=\nabla u_\e$, which corresponds to the situation in which both the surfactant measure $\mu_\e=\rho_\e{\mathcal L}^N$ and the approximating phase interface $\{\nabla u_\e\simeq\frac{1}{\e}\}$ are concentrating on the same $(N-1)$-dimensional set. As explained in \cite{fms} the scaling factor $\e$ multiplying the third integral is chosen in order to observe the effect of this concentration in the asymptotic limit energy. In fact in \cite{fms} the authors have proven that, carrying out the $\Gamma$-limit of ${\mathcal E}_\e$ with respect to the strong $L^1$ convergence of the phase fields and the weak$\ast$-convergence of the surfactant measures, one obtains a limit functional finite for $u\in BV(\Omega, \{ 0,1 \})$ and $\mu\in{\mathcal M}(\Omega)$ (the space of positive Radon measures) where it takes the form
\begin{align}\label{introlimit}
	{\mathcal E}(u,\mu) := \int_{S_{u}\cap\Omega} \Psi \left( \frac{d \mu}{d \HN-1 \llcorner S_{u}} \right)\, d \HN-1 . 
\end{align}
Here $\Psi:[0,+\infty)\to[0,+\infty)$ is a decreasing function of the relative density of the surfactant measure with respect to the surface measure of the interface. In other words, according to the limit energy functional, the surface tension between the phases $\{u=0\}$ and $\{u=1\}$ can be lowered increasing the surfactant density on the interface, a phenomenon that characterizes surface active agents as already recalled at the beginning of this section. It is worth mentioning that in \cite{abcs, acs} such a limit energy has been obtained via a variational discrete-to-continuum coarse-graining procedure starting from the microscopic Blume-Emery-Griffiths ternary surfactant model.\\

In this paper we are interested in extending the results above to the framework of solid-solid phase transition models. Unlike the fluid-fluid ones already introduced, these transitions, and the variational energy models leading to the associated phase separation phenomena, are vectorial problems. The energy functionals we are interested in will be obtained by adding a surfactant-solid interaction term to the functionals $H_\e: W^{2,2}(\Omega;\real^d)$ defined as 
\begin{equation}
	H_\varepsilon(u) := \int \limits_\Omega \frac{1}{\varepsilon}W(\nabla u)
 	+ \varepsilon |\nabla^2 u|^2 \,dx.
\end{equation}
The latter functionals, which are the analogue for gradient vector fields of the Cahn-Hillard functionals $CH_\e$ mentioned before, commonly arise as higher-order regularizations of non-convex stored energy functional in elasticity as those considered in the seminal paper \cite{bj}. Their $\Gamma$-convergence as $\e\to 0$ has been carried out in \cite{cfl} assuming the wells of $W$ to be rank-$1$ connected; i.e., 
\begin{equation} \label{tempKorr5}
\{W=0\}=\{A,B\}, \text{ with } A-B=a\otimes\nu, \text{ for some }a\in\real^{d},\,\nu\in{\mathbb S}^{N-1}.
\end{equation}
For further generalizations allowing for frame invariant potentials $W$ see \cite{cs1,cs2,df} (the reader interested in vector valued singularly perturbed problems with higher order gradients regularizations would also find interesting the results obtained in \cite{fm, dmfl, csz, cdmfl}). Under additional assumptions on $W$ and $\Omega$, satisfied in particular by prototypical quadratic potentials as $W(\xi)=\min\{|\xi-A|^2,\,|\xi-B|^2\}$ and by regular convex domains $\Omega$, the authors of \cite{cfl} compute the $\Gamma$-limit of $H_\e$ and prove that the latter is given by a functional $H$ finite on those $u\in W^{1,1}(\Omega;\real^d)$ with $\nabla u\in BV(\Omega;\{A,B\})$. On this set of functions $H(u)$ takes the form
\begin{equation}
H(u)=K\cdot \HN-1(S_{\nabla u}\cap\Omega),
\end{equation}
where the constant $K>0$ is obtained by solving an asymptotic cell-problem formula.\\ 

In the present paper we are going to investigate functionals defined on functions $u \in W^{2,2}(\Omega, \real^d)$ and measures $\mu = \rho \mathcal{L}^N$ of the form 
\begin{equation} \label{tempKorr4}
	E_\varepsilon(u,\mu) := \int \limits_\Omega \frac{1}{\varepsilon}W(\nabla u)
 	+ \varepsilon |\nabla^2 u|^2 + \varepsilon (\rho - |\nabla^2 u|)^2 \,dx.
\end{equation}
The main result of this paper is stated in  Theorem \ref{mainthm} in which we compute, under the same assumptions on $W$ (w.l.o.g. we assume that in \eqref{tempKorr5} $\nu=e_N$) and $\Omega$ as those considered in \cite{cfl} (see Section \ref{sec:enerfunc} for details), the $\Gamma$-limit of $E_\varepsilon(u,\mu)$ with respect to the strong $W^{1,p}$ convergence of the deformations $u$ and the weak$\ast$-convergence of the surfactant measures $\mu$. We have that $\gammalim_{\varepsilon \rightarrow 0^+} E_\varepsilon(u,\mu) = E(u,\mu)$ where $E(u,\mu)$ is a functional finite on those functions $u\in W^{1,p}(\Omega;\real^d)$ such that $\nabla u\in (\Omega;\{A,B\})$ on which it takes the form 
\begin{equation}\label{intro:energy}
	E(u,\mu) = \int_{S_{\nabla u}\cap\Omega} \Phi \left( \frac{d \mu}{d \HN-1 \llcorner S_{\nabla u}} \right)\,d \HN-1.
\end{equation}
The surface tension $\Phi$ above is a nonnegative nonincreasing function given by an asymptotic formula (see Definition \ref{definitionPhi}). 
Roughly speaking $\Phi(\gamma)$ can be interpreted as
the effective energy per unit $\HN-1$-measure associated to the phase separation induced by the deformation $u: Q \rightarrow \real^d$ with
\begin{align*}
	\nabla u(x) = 
	\begin{cases} 
		-a \otimes e_N ,& x_N < 0,\\
	a \otimes e_N, &x_N > 0
	\end{cases}
\end{align*}
and in presence of the surfactant measure
\begin{equation*}
	\mu := \gamma \HN-1 \llcorner \{x_N = 0\}.
\end{equation*}
In order to prove our main result, we first need to show that $\Phi$ can be obtained restricting the class of admissible functions in the asymptotic formula to those sharing additional regularity and periodicity assumption as in \eqref{Phiptilde}. To this end we need to combine some of the arguments in \cite{cfl} and \cite{fms}, the latter modified to fit in the present vectorial case. 
Such characterization allows us to compute the $\Gamma$-limit on functions with
fixed boundary conditions that are periodic in direction of the phase separation. The proof of the $\gammaliminf$-inequality (Proposition \ref{gammaLimInf}) is then obtained by a blow-up technique near the interfaces of $u$. The proof of the $\Gamma$-limsup inequality (Proposition \ref{recoverySequence}) makes use of a density argument which reduces the construction of a recovery sequence for a generic pair $(u,\mu)$ to the case of deformations with a single interface and to constant surfactant densities.

\section{Notation and Preliminaries}

\subsection{General notation}

Let $\Omega \subset \real^N$ be a bounded, open set with Lipschitz boundary. We denote by ${\mathcal L}
^N$ and ${\mathcal H}^{N-1}$ the $N$-dimensional Lebesgue measure and the $(N-1)$-dimensional
 Hausdorff measure in $\mathbb R^N$, respectively. We use the notation 
$|U| := \int_U \,dx := \mathcal{L}^N(U)$.  We denote by ${\mathcal M}(\Omega)$ the space of 
non negative Radon measures finite on $\Omega$. 
We set $Q_{x_0,\delta} := (x_0 - \delta, x_0 + \delta)^N$ and we use the notation $Q := Q_{0,\frac{1}{2}}$
for the unitary open cube in $\R^N$ centred at the origin. Given $x\in\real^N$ we label the first 
$(N-1)$-coordinates as 
$x'$ and the last coordinate as $x_N$ and we write $x=(x',x_N)$. 
We also set $Q' = \big( - \frac{1}{2}, \frac{1}{2} \big)^{N-1}$, 
hence $Q = Q' \times \big( - \frac{1}{2}, \frac{1}{2} \big)$.
Given a function $u: Q \rightarrow \real^d$ such that for all $x_N \in \big( - \frac{1}{2}, \frac{1}{2} \big)$
it holds that $u(\cdot,x_N)$ is $Q'$-periodic, we say that $u$ is $Q'$-periodic in $x'$.
Given a function $u \in L^1(\Omega, \real^d)$, we denote by $S_u$ the approximate discontinuity set of $u$,
i.e., the set of those points $x \in \Omega$ for which no $z \in \real^d$ exists such that
$\lim_{r \rightarrow 0^+} |B_r(x)|^{-1} \int_{B_r(x)} |u(y) -z| dy = 0$ holds. We denote by $BV(\Omega)$ 
the set of functions of bounded variation in $\Omega$. 
We say that a measurable set $E\subset \R^N$ is a set of finite perimeter in
 $\Omega$ if $\chi_E\in BV(\Omega)$. Denoting by $P(E,\Omega)$ the De Giorgi's perimeter in $\Omega$ of
 $E$, if $E$ is a set of finite perimeter we also write that
$P(E,\Omega)={\mathcal H}^{N-1}(\partial^*E\cap\Omega)<+\infty$ where $\partial^*E$ 
stands for the reduced boundary of $E$.
If $u\in BV(\Omega;\{a,b\})$ is a function of bounded variation in $\Omega$ taking only the 
two values $a,b\in\R^d$, the $(N-1)$-Hausdorff measure
of $S_u$ equals the perimeter of the level set $\{u = a\}$ (and $\{u = b\}$) in $\Omega$ or in formula 
$\HN-1(S_u) = P(\{u = a\}, \Omega)$. For all properties of functions of bounded variations and of 
sets of finite perimeter needed in this paper we refer the reader to \cite{afp}. Finally we set 
${\mathbb S}^{N-1}:=\{x\in\real^N:\,|x|=1\}$ and we denote by $c$ and $C$ generic real
 positive constants that may vary from line to line and expression to expression within the same formula.

\subsection{Preliminaries}
In what follows we will often make use of the following straightforward identity: 
\begin{equation}\label{(3.5)}
\min \{\lambda^2 + w^2, 2w^2\} = w^2 +(\max \{\lambda + w,0\} - w)^2
\text{ for all } \lambda \leq 0, w \geq 0. \end{equation}
The next theorem is proved in \cite[Theorem 3.3]{cfl}. 
\begin{thm}\label{laminates}
Let $u \in W^{1,1}(\Omega, \real^d)$ be such that $\nabla u \in \BV(\Omega,\{A,B\})$ with
 $A-B=a\otimes\nu, \text{ for some }a\in\real^{d},\,\nu\in{\mathbb S}^{N-1}$.
Then $u$ has the form
\begin{equation*}
	u(x',x_N) = \gamma_0 + ax_N - 2 \psi(x)a
\end{equation*}
for some $\gamma_0 \in \real^d$ such that $\gamma_0 \cdot a = 0$ and for some
$\psi \in W^{1,\infty}(\Omega, \real^d)$ such that $\nabla \psi(x) = \chi_E(x) e_N$. The set
$E \subseteq \Omega$ has $P(E,\Omega) < \infty$ and moreover
\[\partial^\ast E = \bigcup_{i=1}^\infty \Omega_i \times \{t_i\}\]
with $\Omega_i \subset \real^{N-1}$ connected, open and bounded and
$t_i \in \real$.
If in addition for each $t \in \real$ the set
$\{(x',x_N) \in \Omega \mid x_N=t\}$ is connected, then $\psi$ depends
only on $x_N$.
\end{thm}

The following lemma is proved in \cite[Lemma 3.2]{fms}. 
\begin{lem}\label{L3.2}
Assume that $(X,\mu)$ is a measure space with $\mu$ a non-atomic positive
measure. Let $g:X \rightarrow [0,\infty) \in L^1(X,\mu) \cap L^2(X,\mu)$ and
$0 < \gamma \leq \int_X g \,d\mu$ be given.
Then, for all $v \geq 0$ with $\int_X v\,d\mu = \gamma$ it holds true that
\begin{equation*}
	\int_X (v-g)^2 \,d\mu \geq \int_X (\max \{ \lambda + g, 0 \} - g)^2 \,d\mu 
		= \int_X \min\{ \lambda^2, g^2 \} \,d\mu,
\end{equation*}
where $\lambda \in (-\infty,0]$ satisfies $\int_X \max\{\lambda + g,0\}\, d\mu = \gamma$.
\end{lem}

\subsection {The energy functional}\label{sec:enerfunc}
In this section we introduce the energy functional we are interested in. For $\varepsilon>0$ we consider the functional 
$
E_\varepsilon:W^{1,1}(\Omega;\real^d)\times\mathcal M(\Omega)\times{\mathcal A}(\Omega)\to[0,+\infty]
$
defined as

\begin{align}
	E_\varepsilon(u,\mu,U) := 
	\begin{cases} 
		\int \limits_U \frac{1}{\varepsilon}W(\nabla u)
 		+ \varepsilon |\nabla^2 u|^2 + \varepsilon (\rho - |\nabla^2 u|)^2 \,dx, &\text{if }
		u\in W^{2,2}(\Omega;\real^d), \,\mu=\rho\, \,dx,\\
	+\infty, &\text{otherwise in }W^{1,1}(\Omega;\real^d)\times\mathcal M(\Omega)
	\end{cases}
\end{align}
With a little abuse of notation we will also introduce the functional 
$E_\varepsilon(u,\mu):W^{1,1}(\Omega;\real^d)\times\mathcal M(\Omega)\to [0,+\infty]$ defined as  
\begin{equation*}
	E_\varepsilon(u,\mu):=E_\varepsilon(u,\mu,\Omega).
\end{equation*}
The asymptotic analysis as $\varepsilon\to 0$ of the functional $E_\varepsilon(u,\mu)$ will be carried
over in the ambiance space $W^{1,1}(\Omega)\times\mathcal M(\Omega)$ endowed with the convergence 
$\tau_1\times\tau_2$ where $\tau_1$ denotes the strong convergence in $W^{1,1}(\Omega;\real^d)$,
 while $\tau_2$ denotes the weak*-convergence in the space of non-negative bounded 
Radon measures ${\mathcal M}(\Omega)$. \\

On the potential $W:\real^{d\times N} \rightarrow [0,\infty)$ we make the following set of assumptions:
\begin{align*} \tag{H1} \label{(H1)}
	&W \text{ is continuous, } W(\xi)=0  \text{ if and only if }
		\xi\in\{A ,B\} \\
	&\text{ where } A-B=a\otimes\nu, \text{ for some }a\in\real^{d},\,\nu\in{\mathbb S}^{N-1}.
\end{align*}
\begin{equation*} \tag{H2} \label{(H2)''} 
	\frac{1}{C} |\xi|^p - C \leq W(\xi) \leq C(|\xi|^p +1) 
	\text{ for some } C > 1,  p\geq 2
\end{equation*}
\begin{align*} \tag{H3} \label{(H4)}
	&c |\xi - A|^p \leq W(\xi) \leq C|\xi - A|^p,
		\quad |\xi - A| \leq \rho \nonumber \\
	&c |\xi - B|^p \leq W(\xi) \leq C|\xi - B|^p,
		\quad |\xi - B| \leq \rho \nonumber \\
	&\text{for some } \rho > 0 \text{ and } p\geq 2
\end{align*}
\begin{equation*} \tag{H4} \label{(H5)}
	W(\xi_1,\dots, \xi_i, \dots, \xi_N) = W(\xi_1,\dots, -\xi_i, \dots, \xi_N),
	\quad i= 1, \dots, N
\end{equation*}

\begin{rem}
We observe that assumption (H1) together with the control from below in (H2) for $p>1$ would suffice to obtain the forthcoming compactness and liminf inequality statements.
\end{rem}

\begin{rem} \label{R6.1}
The following control on the potential energy $W$ is proven in \cite[Remark 6.1]{cfl}. From \eqref{(H2)''} and \eqref{(H4)}, there exist $C_1,C_2 > 0$ such that
\begin{equation*}
	C_1 |\xi'|^p \leq W(\xi) \leq C_2 (W(\eta) + |\xi - \eta|^p) 
\end{equation*}
holds for all $\xi,\eta \in \real^{d \times N}$.
\end{rem}
\section{Compactness and characterization of the asymptotic surface tension}
In this section we prove the compactness statement for sequences $u_h$ and $\mu_h=\rho_h\,dx$ with equibounded energy $E_{\e_h}(u_h.\rho_h)$. We moreover introduce the effective asymptotic (as $\e_h\to 0$) surface tension of the energy functional given by the variational limit of $E_{\e_h}$ and provide some useful characterization of it. For simplicity of notation and without loss of generality  from now on we will assume in $(H1)$ that $A = -B = a \otimes e_N$.\\

\subsection{Compactness}
In what follows we state the main compactness theorem for our functionals. In all our analysis $\e_h$ denotes a sequence of positive numbers vanishing as $h\to +\infty$.

\begin{thm}[Compactness]\label{compactness}
Let $W:\real^{d\times N} \rightarrow [0,\infty)$ satisfy assumptions \eqref{(H1)} and \eqref{(H2)''}, $\varepsilon_h \rightarrow 0^+$,
$(u_h)$ be a sequence in $W^{2,2}(\Omega, \real^d)$ and $\mu_h = \rho_h \,dx$ be a
sequence in $\mathcal{M}(\Omega)$ such that 
\begin{equation*}
	\sup_h E_{\varepsilon_h}(u_h,\rho_h) < \infty \text{ and }
		\sup_h \mu_h(\Omega) < \infty.
\end{equation*}
Then there exist subsequences $(u_{h_k})$, $(\rho_{h_k})$ and
$u \in W^{1,p}(\Omega,\real^d)$ with $\nabla u \in \BV(\Omega,\{A,B\})$
and $\mu \in \mathcal{M}(\Omega)$ such that
\begin{equation*}
	 u_{h_k} - \frac{1}{|\Omega|} \int_\Omega u_{h_k} \,dx \rightarrow u 
	\text{ in } W^{1,p}(\Omega, \real^d) \text{ and } \mu_{h_k} \weak* \mu \text{ as }k\to +\infty.
\end{equation*}
\end{thm}
\begin{proof}
The convergence of a subsequence of $(u_h)$ follows as in \cite[Remark 3.2 (ii)]{cfl}. In fact we can apply the proof of \cite[Theorem 3.1]{cfl} to find a subsequence $u_{h_k} \rightarrow u$ in $W^{1,1}(\Omega,\real^d)$.
By assumption (\ref{(H2)''}) there exists $L > 0$ such that $W(\xi) \geq c |\xi|^p$ for $p\geq 2$ and for all $\xi$ satisfying $|\xi| \geq L$
and therefore
\begin{equation*}
	c\int_{ \{ |\nabla u_{h_k}| \geq L \} } |\nabla u_{h_k}|^p \,dx 
	\leq \int_\Omega W(\nabla u_{h_k}) \,dx \rightarrow 0.
\end{equation*}
If follows
\begin{equation*}
	\int_{ \{ |\nabla u_{h_k}| \geq L \} } |\nabla u_{h_k} - \nabla u|^p \,dx \rightarrow 0
\end{equation*}
and together with
\begin{equation*}
	\int_{ \{ |\nabla u_{h_k}| \leq L \} } |\nabla u_{h_k} - \nabla u|^p \,dx
	\leq	 (L + |A| + |B|)^{p-1} \int_\Omega |\nabla u_{h_k} - \nabla u| \,dx \rightarrow 0
\end{equation*}
it implies the convergence of $u_{h_k}$ in $W^{1,p}(\Omega,\real^d)$.
The convergence of a subsequence of $(\mu_h)$ in the weak$^\ast$-topology is a consequence of the weak$^\ast$-compactness 
of $\mathcal{M}(\Omega)$.
\end{proof}

In what follows we define the lower semicontinuous envelope of our energy functional on a restricted class of admissible functions. In Section \ref{gammalim} we are going to prove that this lower bound is actually our $\Gamma$-liminf functional.\\

We start by introducing a notation. For every $U \subset \Omega$ open set we define
\begin{equation*}
	F_\varepsilon(u,\lambda, U) := 
	\int_U\frac{1}{\varepsilon} W(\nabla u) + \varepsilon |\nabla^2 u|^2 
	+ \varepsilon \min \{\lambda^2, |\nabla^2 u|^2 \}\,\,dx.
\end{equation*}
In the case $U=Q$ (the unitary cube in $\R^N$ centred at the origin) the notation above will be shortened and we will use $F_\varepsilon(u,\lambda) := F_\varepsilon(u,\lambda, Q)$. 

\begin{df}\label{definitionPhi}
For $\gamma \geq 0$, $k > 0$ and 
$\omega \subset \real^{N-1}$ a bounded open set with ${\mathcal H}^{N-1} (\partial \omega) = 0$  we define
\begin{align*}
		F(\gamma, \omega \times (-k,k)) := \inf \Big\{ \liminf_{h \rightarrow \infty} 
		F_{\varepsilon_h}(u_h,\lambda_h, \omega \times (-k,k)): \,
		u_h \rightarrow |x_N|a, (u_h, \lambda_h) \in 
		\mathcal{A} (\gamma, \omega \times (-k,k)) \Big \}
\end{align*}
where we have used the notation $\mathcal{A}(\gamma, \omega \times (-k,k))$ for the set of admissible functions
defined as
\begin{equation*}
	\mathcal{A}(\gamma, \omega \times (-k,k)) := \Big\{ (u, \lambda) 
	\in W^{2,2}(Q,\real^d) \times (-\infty, 0] \mid 
	\int_{\omega \times (-k,k)} \max \{ \lambda + |\nabla^2 u|,0 \} \,\,dx
	\leq \gamma\, {\mathcal H}^{N-1}(\omega) \Big\}
\end{equation*}
\end{df}

In order to shorten the notation we also set
\begin{equation} \label{tempKorr6}
\mathcal{A}(\gamma) := \mathcal{A}(\gamma,Q), \quad
	\tilde \Phi(\gamma) := F(\gamma,Q) \quad \text{and} \quad
	\Phi(\gamma) := \lim_{\delta  \rightarrow 0^+} \tilde \Phi(\gamma + \delta)
\end{equation}
and we observe that since $\tilde \Phi$ is non-increasing, the function $\Phi$
is well-defined.\\

The following lemma, whose proof we omit, can be proved as in \cite[Lemma 4.3]{cfl}, the only care being that the rescaling argument used to show assertion (iv) now makes use of the admissible sequence $(\alpha u_n(\frac{x}{\alpha}), \frac{\lambda_n}{\alpha})$.
\begin{lem}\label{L4.3}
For $\gamma \geq 0$ fixed, it holds
\begin{align*}
	(i) &\quad F(\gamma, x' + \omega \times (-k,k)) = F(\gamma, \omega \times (-k,k))
		\text{ for all } x' \in \real^{N-1} \\
	(ii)& \quad  F(\gamma, \omega_1 \times (-k,k)) \leq F(\gamma, \omega_2 \times (-k,k))
		\text{ if } \omega_1 \subset \omega_2 \\
	(iii)& \quad F(\gamma, \omega_1 \times (-k,k)) + F(\gamma, \omega_2 \times (-k,k))
		\leq F(\gamma, \omega_1 \cup \omega_2 \times (-k,k)) \text{ if }
		\omega_1 \cap \omega_2 = \emptyset \\
	(iv)& \quad F(\gamma, \alpha \omega \times (- \alpha k, \alpha k)) 
		=\alpha^{N-1} F(\gamma, \omega \times (-k,k)) \text{ for } \alpha > 0, \\
	&\quad F(\gamma, \alpha \omega \times (-h, h)) 
		\geq \alpha^{N-1} F(\gamma, \omega \times (-k,k)) \text{ for } 1 > \alpha > 0 \\
	(v)& \quad F(\gamma, \omega \times (-k,k)) = {\mathcal H}^{N-1}(\omega) F(\gamma, Q' \times (-k,k)) \\
	(vi)& \quad F(\gamma, \omega \times (-k,k)) = F(\gamma, \omega \times (-k',k'))
		\text{ for all } k' > 0
\end{align*}
and analogously for $\lim_{\delta \rightarrow 0^+} F(\gamma + \delta, \omega \times (-k,k))$.\\
In particular it holds that $F(\gamma, \omega \times (-k,k)) = {\mathcal H}^{N-1} (\omega) \tilde \Phi(\gamma)$.
\end{lem}

\subsection {Characterization of the surface tension}
In this section we further characterize the surface tension $\Phi$ and $\tilde \Phi$. More specifically, following \cite{cfl}, we prove that the minimum problem in Definition \ref{definitionPhi} can be restrict to a narrower class of competitors. 

\begin{prop} \label{verticalMatching}
Let $W:\real^{d\times N} \rightarrow [0,\infty)$ satisfy \eqref{(H1)}, \eqref{(H2)''} 
and \eqref{(H4)} and let
$\gamma \geq 0$ and $\delta > 0$ be given. Then there exist sequences
$\varepsilon_h^\delta \rightarrow 0^+$ and $(u_h^\delta, \lambda_h^\delta)
\in \mathcal{A} (\gamma + \delta)$ satisfying
\begin{align*}
	& u_h^\delta \rightarrow |x_N|a\, \text{ in } W^{1,p} (Q,\real^d), \\
	& u_h^\delta = -x_N a \,\text{ near } x_N = -\frac{1}{2}, \quad
		u_h^\delta = x_N a + c_h^\delta \text{ near } x_N = \frac{1}{2},
		\quad c_h^\delta \rightarrow 0 \text{ as } h \rightarrow \infty, \\
	& \lim_{h \rightarrow \infty} F_{\varepsilon_h}
		(u_h^\delta, \lambda_h^\delta) = \tilde \Phi(\gamma).
\end{align*}
\end{prop}
\begin{proof}
Our proof follows the strategy of the proof of \cite[Proposition 6.2]{cfl}.
We fix $\delta$ and drop it from the notation.
Choosing admissible sequences satisfying
\begin{equation}\label{optimalbound}
\lim_{h\rightarrow \infty} F_{\varepsilon_h} (u_h,\lambda_h) = \tilde \Phi(\gamma)
\end{equation}
and using  the compactness result \eqref{compactness}
we can assume that $u_h \rightarrow |x_N|a =: u_0$ in $W^{1,p}(Q,\real^d)$.
We can partition $Q' \times (\frac{1}{6}, \frac{1}{3})$
 into $\lfloor \frac{1}{\varepsilon_h} \rfloor$-horizontal layers
of height $\lfloor \frac{1}{\varepsilon_h} \rfloor^{-1} \frac{1}{6}$ and choose
a layer $L_h = Q' \times (\theta_h - \lfloor \frac{1}{\varepsilon_h} \rfloor^{-1} \frac{1}{6}, \theta_h)$
which satisifies
\begin{align} \label{temp11}
	 \left\lfloor \frac{1}{\varepsilon_h} \right\rfloor &\left( F_{\varepsilon_h}(u_h,\lambda_h,L_h)
		+ \int_{L_h} |\nabla u_h - a \otimes e_N| ^p + |u_h - u_0(x)|^p \,\,dx \right) \nonumber \\
	&\leq F_{\varepsilon_h}\left(u_h,\lambda_h,Q' \times \left(\frac{1}{6}, \frac{1}{3}\right) \right)
		+ \int_{Q' \times (\frac{1}{6},\frac{1}{3})} 
		|\nabla u_h - a \otimes e_N| ^p + |u_h - u_0(x)|^p \,\,dx \nonumber \\
	&	=: \alpha_h \rightarrow 0,
\end{align}
where we have used \eqref{L4.3} (vi) which asserts that the energy concentrates near $Q'\times\{0\}$.
By the continuity of $W$, \eqref{optimalbound} and the very definition of the energy functional
$F_{\varepsilon}$ there exists $z_h \in (\theta_h - \lfloor \frac{1}{\varepsilon_h} \rfloor^{-1}
 \frac{1}{6}, \theta_h)$ such that
\begin{align} \label{temp12}
	\int_{Q'} &\frac{1}{\varepsilon_h} W(\nabla u_h(x',z_h)) 
		+ \varepsilon_h |\nabla^2 u_h(x',z_h)|^2 + \varepsilon_h
		\min \{ \lambda_h^2, |\nabla^2 u_h(x',z_h)|^2 \} \,\,dx'\nonumber \\
	& +\int_{Q'} |\nabla u_h(x',z_h) - a \otimes e_N| ^p 
		+ |u_h(x',z_h) - u_0(x',z_h)|^p\,\,dx' \leq 6 \alpha_h.
\end{align}
Choosing a smooth cut-off function $\varphi_h: \left( -\frac{1}{2},\frac{1}{2}\right) \rightarrow \real$ 
satisfying $\varphi_h(x_N) = 1$ if
$x_N \leq \theta_h - \lfloor \frac{1}{\varepsilon_h} \rfloor^{-1} \frac{1}{6}$, 
$\varphi_h = 0$ if $x_N \geq \theta_h$, 
$|\varphi_h'| \leq \frac{c}{\varepsilon_h}, |\varphi_h''| \leq \frac{c}{\varepsilon_h^2}$,
we define
\begin{align*}
	v_h(x) &:=  u_0(x) + u_h(x',z_h) - u_0(x',z_h) 
		+ \varphi_h(x_N)( u_h(x) - u_0(x) -(u_h(x',z_h) - u_0(x',z_h)) ) \\
	&= \varphi_h(x_N) u_h(x) + (1 - \varphi_h(x_N)) (u_0(x) + u_h(x',z_h) - u_0(x',z_h)).
\end{align*}
We claim that the following limits hold true:
\begin{align*}
	&(a)\int_{L_h} |v_h - u_0|^p \,\,dx \rightarrow 0, &(b)&\, \frac{1}{\varepsilon_h} \int_{L_h}
		|\nabla v_h - a \otimes e_N|^p \,\,dx \rightarrow 0, \\
	&(c)\,F_{\varepsilon_h} (v_h,\lambda_h,L_h) \rightarrow 0, &(d)& \int_{L_h} \max \{ \lambda_h 
		+ |\nabla^2 v_h|, 0 \} \,\,dx \rightarrow 0.
\end{align*}
The limit in $(a)$ follows directly from (\ref{temp11}) and (\ref{temp12}).
We now prove the limit in $(b)$. We can apply Poincar\'e inequality to the function 
$u_h(x) - u_0(x) - (u_h(x',z_h) - u_0(x',z_h))$
to obtain
\begin{align} \label{temp13}
	\frac{1}{\varepsilon_h}& \int_{L_h} |\nabla v_h - a \otimes e_N|^p \,\,dx \nonumber\\
	&\leq \frac{C}{\varepsilon_h} \int_{L_h} |\nabla u_h(x) - a \otimes e_N|^p 
		+ |\nabla_{x'}u_h(x',z_h)|^p+ \frac{c^p}{\varepsilon_h^p}
		 |u_h(x) - u_0(x) - (u_h(x',z_h) - u_0(x',z_h))|^p \,\,dx  \nonumber\\
	&\leq C\int_{Q'} |\nabla_{x'} u_h(x',z_h)|^p \,\,dx' 
		+ \frac{C}{\varepsilon_h} \int_{L_h} |\nabla u_h - a \otimes e_N|^p \,\,dx \nonumber \\
	&\leq C \int_{Q'} |\nabla u_h(x',z_h) - a \otimes e_N|^p \,\,dx' 
		+ \frac{C}{\varepsilon_h} \int_{L_h} |\nabla u_h - a \otimes e_N|^p \,\,dx
		\rightarrow 0,
\end{align}
where in the last step we have used (\ref{temp11}) and (\ref{temp12}). We now prove claim $(c)$.
Thanks to \eqref{(H2)''} and \eqref{(H4)} we obtain
\begin{align} \label{temp14}
	\frac{1}{\varepsilon_h} \int_{L_h} W(\nabla v_h) \,\,dx &\leq \frac{1}{\varepsilon_h}
		\int_{L_h \cap \{ |\nabla v_h - \nabla u_0| < \rho \} } C|\nabla v_h - a \otimes e_N|^p \,\,dx
		+ \frac{1}{\varepsilon_h} \int_{L_h \cap \{ |\nabla v_h - \nabla u_0| \geq \rho \} }
		C(1 + |\nabla v_h|^p) \,\,dx \nonumber \\
	 &\leq \frac{C}{\varepsilon_h} \int_{L_h} |\nabla v_h - a \otimes e_N|^p \,\,dx
		\to 0.
\end{align}
where the last limit follows by $(b)$. We also have that
\begin{align*}
	\int_{L_h}  \varepsilon_h |\nabla^2 v_h|^2 \,\,dx \leq& C \int_{L_h} \varepsilon_h |\nabla_{x'} u_h(x',z_h)|^2
		+ \varepsilon_h |\nabla^2 u_h|^2 + \frac{1}{\varepsilon_h}|\nabla u_h - a \otimes e_N|^2 \\
	& \quad + \frac{1}{\varepsilon_h^3} | u_h - u_0 - (u_h(x',z_h) - u_0(x',z_h))|^2 \,\,dx \\
	\leq& C \varepsilon_h^2 \int_{Q'} |\nabla_{x'}^2 u_h(x',z_h)|^2 \,\,dx' + C\varepsilon_h
		\int_{L_h} |\nabla^2 u_h|^2 \,\,dx \\ &
		+ \frac{C}{\varepsilon_h} \int_{L_h} |\nabla u_h - a \otimes e_N|^2 \,\,dx 
	+ C \int_{Q'} |\nabla_{x'} u_h(x',z_h)|^2 \,\,dx' \\&+ \frac{C \varepsilon_h^{(p-2)/p}}{\varepsilon_h^3}
		\left( \int_{L_h} |u_h(x) - u_0(x) - (u_h(x',z_h) - u_0(x',z_h))|^p \,\,dx \right) ^{2/p} \\
	\leq &C \varepsilon_h^2 \int_{Q'} |\nabla_{x'}^2 u_h(x',z_h)|^2 \,\,dx' + C\varepsilon_h
		\int_{L_h} |\nabla^2 u_h|^2 \,\,dx 
		+\frac{C}{\varepsilon_h^{2/p}} \left( \int_{L_h} |\nabla u_h - a \otimes e_N|^p \,\,dx \right) ^{2/p} \\
	&\quad + C \left( \int_{Q'} |\nabla u_h(x',z_h) - a \otimes e_N|^p \,\,dx' \right) ^{2/p} \rightarrow 0
\end{align*}
where we have used H\"older's inequality, Poincar\'e's inequality for $u_h(x) - u_0(x) - (u_h(x',z_h) - u_0(x',z_h))$
and, in the last step, (\ref{temp11}), (\ref{temp12}) and (\ref{temp13}). By the trivial inequality $ \int_{L_h}
\varepsilon_h \min \{\lambda_h^2, |\nabla^2 v_h|^2 \}\,\,dx\leq \int_{L_h}\e_h |\nabla^2 v_h|^2\,dx$,
the estimate above shows $(c)$.
Finally $(d)$ follows from (\ref{temp11}) and (\ref{temp12}) thanks to the estimate
\begin{equation*}
	\int_{L_h} \max \{ \lambda_h + |\nabla^2 v_h|, 0 \} \,\,dx \leq \int_{L_h} |\nabla^2 v_h| \,\,dx
	\leq C \varepsilon_h^{1/2} \left( \int_{L_h} |\nabla^2 v_h|^2 \,\,dx \right)^{1/2} \rightarrow 0.
\end{equation*}
In the next step, we choose a smooth  cut-off function $\psi: \left( -\frac{1}{2},\frac{1}{2}\right) \to\R$ such that 
$\psi(x_N) = 1$ for $x_N \leq \frac{1}{3}$ and $\psi = 0$  near  $x_N = \frac{1}{2}$. We moreover 
may assume that $\|\psi'\|_\infty\leq C$ and $\|\psi''\|_\infty\leq C$ and define
\begin{equation*}
w_h(x) := u_0(x) + c_h + \psi(x_N) (u_h(x',z_h) - u_0(x',z_h) - c_h)
\end{equation*}
where $c_h := \int_{Q'} u_h(x',z_h) - u_0(x',z_h) \,\,dx' \rightarrow 0$ by (\ref{temp12}).
We write  $\hat Q_h :=Q' \times (\theta_h, 1/2) $ and claim that
\begin{align*}
	&(a')\int_{\hat Q_h} |w_h - u_0|^p \,\,dx \rightarrow 0, &(b')&\, \frac{1}{\varepsilon_h} \int_{\hat Q_h}
		|\nabla w_h - a \otimes e_N|^p \,\,dx \rightarrow 0, \\
	&(c')\,F_{\varepsilon_h}(w_h,\lambda_h,\hat Q_h) \rightarrow 0, &(d')& \int_{\hat Q_h} 
		\max \{ \lambda_h + |\nabla^2 w_h|, 0 \} \,\,dx
		\leq \int_{\hat Q_h} \max \{ \lambda_h + |\nabla^2 u_h|, 0 \} \,\,dx + d_h, \\
	&\, &\, & \text{where } d_h \rightarrow 0 \text{ as } h \rightarrow \infty.
\end{align*}
The first claim follows from (\ref{temp12}).
In order to prove claim $(b')$ we use Poincar\'e's inequality and Remark \ref{R6.1}. We obtain
\begin{align}\label{eq:vanish}
	\frac{1}{\varepsilon_h}  \int_{\hat Q_h} |\nabla w_h - a \otimes e_N|^p \,\,dx
		&\leq \frac{C}{\varepsilon_h} \int_{\hat Q_h} |\nabla_{x'} u_h(x',z_h)|^p 
		+ |u_h(x',z_h) - u_0(x',z_h) - c_h|^p \,dx \\
	&\leq \frac{C}{\varepsilon_h} \int_{Q'} |\nabla_{x'} u_h(x',z_h)|^p \,dx'
		\leq C \int_{Q'} \frac{1}{\varepsilon_h} W(\nabla u_h(x',z_h)) \,dx' \rightarrow 0 \nonumber
\end{align}
where the last limit follows by (\ref{temp12}).
Similarly, it holds that
\begin{align*}
	\int_{\hat Q_h} \varepsilon_h|\nabla^2 w_h|^2 \,dx &\leq C \varepsilon_h
		 \int_{\hat Q_h} |\nabla_{x'}^2 u_h(x',z_h)|^2 + |\nabla_{x'} u_h(x',z_h)|^2 \,dx \\
	&\leq C \varepsilon_h \int_{Q'} |\nabla_{x'}^2 u_h(x',z_h)|^2 \,dx'
		+ C \varepsilon_h \left( \int_{Q'} |\nabla u_h(x',z_h) - a\otimes e_N|^p \,dx' \right)^{2/p}
		\rightarrow 0.
\end{align*}
Again, as shown for the claim $(c)$, the last two estimates, together with $(b')$ and (\ref{temp14}),
give $(c')$. To prove the last claim we use the estimate
\begin{equation*}
	|\nabla^2 w_h(x)| \leq C(|u_h(x',z_h) - u_0(x',z_h) -c_h| 
	+ |\nabla_{x'}u_h(x',z_h)|) + |\nabla_{x'}^2 u_h(x',z_h)|
\end{equation*}
to find that
\begin{align*}
	\int_{\hat Q_h}  \max \{ \lambda_h + |\nabla^2 w_h|, 0 \} \,dx
		&\leq \int_{\hat Q_h} \max \{ \lambda_h+|\nabla_{x'}^2 u_h(x',z_h)|, 0 \} \,dx + C\int_{\hat Q_h} |\nabla_{x'}u_h(x',z_h)|\, \,dx \\
	& \leq \int_{\hat Q_h} \max \{\lambda _h+ |\nabla^2 u_h|, 0 \} \,dx
		+ C \left(\int_{\hat Q_h} |\nabla_{x'} u_h(x',z_h)|^2 \,dx\right)^{\frac12}.
\end{align*}
The second summand tending to zero as in \eqref{eq:vanish}, $(d')$ holds true.
Let us define
\begin{equation*}
	U_h := \left\{
	\begin{array}{lll}
	u_h, &x_N < \theta_h - \lfloor\frac{1}{\varepsilon_h}\rfloor^{-1} \frac{1}{6},\\
	v_h, &\theta_h - \lfloor\frac{1}{\varepsilon_h}\rfloor^{-1} \frac{1}{6} \leq x_N \leq \theta_h\\
	w_h, &x_N > \theta_h.
	\end{array} \right.
\end{equation*}
Our claims show that $U_h \rightarrow u_0$ in $W^{1,p}(Q,\real^d)$ and that we may assume
$(U_h,\lambda_h) \in \mathcal{A}(\gamma + \delta/2)$ for $h$ large enough.
Since by \eqref{L4.3} (vi) we have that
\begin{equation*}
	F_{\varepsilon_h}\left(u_h,\lambda_h,Q' \times \left(\frac{1}{6},\frac{1}{2}\right)\right) \rightarrow 0,
\end{equation*}
from our claims it also follows that $\lim_{h \rightarrow \infty} F_{\varepsilon_h}(U_h,\lambda_h) 
= \tilde \Phi(\gamma)$.
The proof is completed on observing that the construction above can be repeated on $Q' \times (-\frac{1}{2},0)$.
\end{proof}

In the next proposition, which is analogous to \cite[Proposition 6.3]{cfl}, we show how to further modify the construction of the sequence of functions in the previous proposition to enforce periodic boundary condition in the $x'$ variable. 

\begin{prop} \label{transversalPeriodicity}
Let $W:\real^{d\times N} \rightarrow [0,\infty)$ satisfy \eqref{(H1)}, \eqref{(H2)''} , \eqref{(H4)} and \eqref{(H5)} and let
$\gamma \geq 0$ and $\delta > 0$ be given. Then there exist sequences
$\varepsilon_h^\delta \rightarrow 0^+$ and $(u_h^\delta, \lambda_h^\delta)
\in \mathcal{A} (\gamma + \delta)$ satisfying
\begin{align*}
	&u_h^\delta \in W^{2,\infty}(Q,\real^d), \\
	& u_h^\delta \rightarrow |x_N|a \text{ in } L^1 (Q,\real^d), \\
	& \nabla u_h^\delta = \pm a \otimes e_N \text{ near } x_N 
		= \pm \frac{1}{2} \text{ and } u_h^\delta 
		\text{ is } Q'\text{-periodic for all } x_N, \\
	&\lim_{h \rightarrow \infty} F_{\varepsilon_h}
		(u_h^\delta, \lambda_h^\delta) = \tilde \Phi(\gamma).
\end{align*}
\end{prop}
\begin{proof}
Fixing $\delta$ and dropping it from the notation, we are going to show how to proceed as in the proof of Proposition 6.3 in \cite{cfl} and show that there exist sequences
$\varepsilon_h \rightarrow 0^+, v_h \in W^{2,\infty}(\real^N,\real^d), \lambda_h \leq 0$
such that
\begin{align} \label{tempKorr1}
	&(v_h ,\lambda_h) \in \mathcal{A} \Big( \gamma + \delta, 2Q' \times 
		\Big( -\frac{1}{2},\frac{1}{2}\Big) \Big), \nonumber \\
	&v_h(\cdot,x_N) \text{ is } 2Q' \text{-periodic for all } x_N, \quad
		\nabla v_h = \pm a \otimes e_N \text{ near } x_N = \pm \frac{1}{2}, \nonumber \\
	&\lim_{h \rightarrow \infty} F_{\varepsilon_h} \Big(v_h,\lambda_h, 2Q' \times 
		\Big( -\frac{1}{2},\frac{1}{2}\Big)\Big) = 2^{N-1} \tilde \Phi(\gamma), \nonumber \\
	&\lim_{h \rightarrow \infty} \int_{2Q' \times (-\frac{1}{2},\frac{1}{2})} |v_h - |x_N|a| \,dx = 0.
\end{align}
With (\ref{tempKorr1}) at hand we can extend $v_h$ linearly to $2Q$, define $u_h(x) := \frac{1}{2} v_h(2x)$ and complete the proof noting that
\begin{align*}
	&(u_h ,2\lambda_h) \in \mathcal{A}(\gamma + \delta), \\
	&u_h(\cdot,x_N) \text{ is } Q' \text{-periodic for all } x_N, \quad
		\nabla u_h = \pm a \otimes e_N \text{ near } x_N = \pm \frac{1}{2}, \\
	&\lim_{h \rightarrow \infty} F_{\frac{\varepsilon_h}{2}} (u_h,\lambda_h) = \tilde \Phi(\gamma), \\
	&\lim_{h \rightarrow \infty} \int_Q |u_h - |x_N|a| \,dx = 0.
\end{align*}
In what follows we prove (\ref{tempKorr1}).
\begin{description}
\item[Step 1.] 
We assume that $N = 2$ and set $u_0 := |x_2|a$.
By  \eqref{verticalMatching} we find sequences
$\varepsilon_h \rightarrow 0^+, (u_h, \lambda_h) \in \mathcal{A}(\gamma + \delta/2)$
such that $u_h \rightarrow u_0$ in $L^1(Q,\real^d)$,
$\nabla u_h = \pm a \otimes e_2$ near $x_2 = \pm \frac{1}{2}$ and
$\lim_{h \rightarrow \infty} F_{\varepsilon_h} (u_h,\lambda_h) = \tilde \Phi(\gamma)$.
Moreover, we may assume $u_h \in C^2(Q,\real^d)$: This follows from 
\eqref{(H2)''} and the fact that
\begin{equation*}
	\int_\Omega \max \{ \lambda + |\rho_\varepsilon \ast u|, 0 \} \,dx \rightarrow
	\int_\Omega \max \{ \lambda + |u|, 0 \} \,dx
\end{equation*}
as $\varepsilon \rightarrow 0$ for $u \in L^1(\real^m, \real^n)$ which
can be shown by an application of the Vitali dominated convergence theorem (see \cite{afp} Exercise 1.18). We can therefore conclude
 as in the proof of \cite[Proposition 6.2]{cfl}.
Setting $I_m := ((-\frac{1}{2}, -\frac{1}{2} + \frac{1}{m}) \cup 
(\frac{1}{2} - \frac{1}{m}, \frac{1}{2})) \times (-\frac{1}{2},\frac{1}{2})$,
we have
\begin{equation*}
	\tilde \Phi(\gamma) = \lim_{h \rightarrow \infty} F_{\varepsilon_h}(u_h,\lambda_h)
	\geq \liminf_{h \rightarrow \infty} F_{\varepsilon_h}(u_h,\lambda_h, Q \backslash I_m)
	\geq \tilde \Phi(\gamma) \Big(1 - \frac{2}{m}\Big),
\end{equation*}
and therefore
\begin{equation*}
	\limsup_{h \rightarrow \infty} F_{\varepsilon_h}(u_h,\lambda_h, I_m) \leq  \tilde \Phi(\gamma) \frac{2}{m}.
\end{equation*}
By the compactness result stated in Theorem \ref{compactness} we have $u_h \rightarrow u_0$ in $W^{1,p}(Q,\real^d)$.
Therefore,  for $h$ sufficiently large, it holds that
\begin{equation} \label{temp15}
	F_{\varepsilon_h}(u_h,\lambda_h,I_m) + \int_{I_m} m^p|\nabla u_h - \nabla u_0|^p + |u_h-u_0| \,dx
	\leq \tilde \Phi(\gamma) \frac{3}{m}.
\end {equation}Let us subdivide 
$\big(-\frac{1}{2}, -\frac{1}{2} + \frac{1}{m}\big) \times \big(-\frac{1}{2},\frac{1}{2}\big)$
and $\big(\frac{1}{2} - \frac{1}{m}, \frac{1}{2}\big) \times \big(-\frac{1}{2},\frac{1}{2}\big)$
into $k$ strips of equal width and order them in pairs.
By (\ref{temp15}) we find a pair of strips
\begin{equation*}
	R_{h,m,k}^+ = (b_{h,m,k}, c_{h,m,k}) \times \left(-\frac{1}{2},\frac{1}{2}\right), \quad
	R_{h,m,k}^- = (-b_{h,m,k}, -c_{h,m,k})\times \left(-\frac{1}{2},\frac{1}{2}\right)
\end{equation*}
such that for $h$ sufficiently large it holds that
\begin{equation} \label{temp16}
	F_{\varepsilon_h}(u_h,\lambda_h,R_{h,m,k}^- \cup R_{h,m,k}^+)
	+ \int_{R_{h,m,k}^-  \cup R_{h,m,k}^+} m^p |\nabla u_h - \nabla u_0|^p + |u_h - u_0| \,dx
	\leq \tilde \Phi(\gamma) \frac{3}{mk}.
\end{equation}
In particular we have that, setting $J_{h,m,k}:=(\frac{b_{h,m,k} + c_{h,m,k}}{2},c_{h,m,k})$
\begin{align*}
	&\int\limits_{J_{h,m,k}} \int_{-\frac{1}{2}}^{\frac{1}{2}}
		\frac{1}{\varepsilon_h} W(\nabla u_h(x)) + \varepsilon_h |\nabla^2 u_h(x)|^2
		+ m^p |\nabla u_h(x) - \nabla u_0(x)|^p + |u_h(x) - u_0(x)| \\& \quad
	+\frac{1}{\varepsilon_h} W(\nabla u_h(-x)) + \varepsilon_h |\nabla^2 u_h(-x)|^2
+ m^p |\nabla u_h(-x) - \nabla u_0(-x)|^p + |u_h(-x) - u_0(-x)| \,dx \\
	&\leq \tilde \Phi(\gamma) \frac{3}{mk} 
\end{align*}
and this shows that there exists $a_{h,m,k} \in \left(b_{h,m,k}, \frac{b_{h,m,k} + c_{h,m,k}}{2}\right)$ 
satisfying
\begin{align} \label{temp17}
&	\int_{-\frac{1}{2}}^{\frac{1}{2}} 
		\frac{1}{\varepsilon_h} W(\nabla u_h(a_{h,m,k},x_2)) 
		+ \varepsilon_h |\nabla^2 u_h(a_{h,m,k},x_2)|^2 \nonumber \\
	&\quad\quad\quad+ m^p |\nabla u_h(a_{h,m,k},x_2) - \nabla u_0(a_{h,m,k},x_2)|^p 
		+ |u_h(a_{h,m,k},x_2) - u_0(a_{h,m,k},x_2)|  \nonumber \\
	&\quad\quad\quad+ \frac{1}{\varepsilon_h} W(\nabla u_h(-a_{h,m,k},x_2)) 
		+ \varepsilon_h |\nabla^2 u_h(-a_{h,m,k},x_2)|^2 \nonumber \\
	&\quad\quad\quad+ m^p |\nabla u_h(-a_{h,m,k},x_2) - \nabla u_0(-a_{h,m,k},x_2)|^p 
		+ |u_h(-a_{h,m,k},x_2) - u_0(-a_{h,m,k},x_2)|  \,dx \nonumber \\
	& \quad\leq 6 \tilde \Phi(\gamma).
\end{align}

\begin{figure}
\begin{tikzpicture}[scale=0.85]
	\filldraw[fill=black!5] (2,-3.5) rectangle (3,3.5);
	\filldraw[fill=black!5] (-2,-3.5) rectangle (-3,3.5);	
	\draw [->](-4,0) -- (4,0) node[anchor=west]{$x_1$};
	\draw [->](0,-4) -- (0,4) node[anchor=south]{$x_2$};
	\draw (-3.5,-3.5) rectangle (3.5,3.5);
	\draw (-1,2) node{$Q$};
	\draw (-0.2,3.5) node[anchor=north]{$\frac{1}{2}$};
	\draw (3.7,0.8) node[anchor=north]{$\frac{1}{2}$};
	\draw (1.5,-3.5) -- (1.5,3.5);
	\draw (0.9,1)node[anchor=north]{$\frac{1}{2} - \frac{1}{m}$};
	\draw (-2.7,-3.5) -- (-2.7,3.5);
	\draw (2.7,-3.5) -- (2.7,3.5);
	\draw (0.7,-0.5)node[anchor=north]{$b_{h,m,k}$};
	\draw [->] (0.7,-0.6)  -- (2,0);
	\draw (0.7,-1.2)node[anchor=north]{$a_{h,m,k}$};
	\draw [->] (0.7,-1.3)  -- (2.7,0);
	\draw (0.7,-2)node[anchor=north]{$c_{h,m,k}$};
	\draw [->] (0.7,-2.1)  -- (3,0);
	\draw (2,3.5) -- (2,3.8);
	\draw (3,3.5) -- (3,3.8);
	\draw[<->] (2,3.7) -- (3,3.7);
	\draw (2.6,4) node{$\frac{1}{mk}$};
\end{tikzpicture}
\caption{Sketch of the construction in Step 1}
\end{figure}
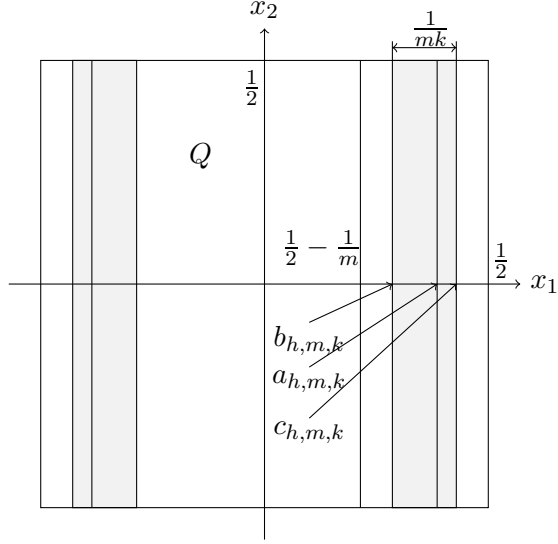

Next we modify $u_h$ in order to obtain a new function that coincides with
$u_h$ near $x_1 = -a_{h,m,k}  + \frac{1}{2mk}$
and with $u_h(-a_{h,m,k},\cdot)$ near $x_1 = -a_{h,m,k}$
(note that by construction $(-a_{h,m,k}, -a_{h,m,k} +  \frac{1}{2mk}) \subset (-c_{h,m,k}, -b_{h,m,k})$).
To this end we choose a smooth cut-off function
$\varphi_{h,m,k}: \left( -\frac{1}{2},\frac{1}{2}\right) \rightarrow \real$ such that
$\varphi_{h,m,k} = 1$ if $x_1 > -a_{h,m,k} + \frac{1}{2mk}$,
$\varphi_{h,m,k} = 0$ if $x_1 < -a_{h,m,k} $,
\begin{equation*}
	|\varphi_{h,m,k}'|_\infty \leq cmk, \quad
	|\varphi_{h,m,k}''|_\infty \leq cm^2k^2
\end{equation*}
and define
\begin{equation*}
w_{h,m,k}(x) :=\varphi_{h,m,k}(x_1) u_h(x) + (1 - \varphi_{h,m,k}(x_1)) u_h(-a_{h,m,k},x_2).
\end{equation*}
We have that $w_{h,m,k} \in W^{2,\infty} (Q, \real^d)$ and $\nabla w_{h,m,k} = \pm a \otimes e_2$
near $x_2 = \pm \frac{1}{2}$. We are going to show that
\begin{align} \label {temp18}
	&\limsup_{h \rightarrow \infty} \limsup_{k \rightarrow \infty} \limsup_{m \rightarrow \infty}
		F_{\varepsilon_h}(w_{h,m,k}, \lambda_h) \leq \tilde \Phi(\gamma), \nonumber \\
	&\limsup_{h \rightarrow \infty} \limsup_{k \rightarrow \infty} \limsup_{m \rightarrow \infty}
	\int_Q |w_{h,m,k} - u_0| \,dx = 0 \quad \text{and} \nonumber \\
	&\limsup_{h \rightarrow \infty} \limsup_{k \rightarrow \infty} \limsup_{m \rightarrow \infty}
	\int_Q \max \{ \lambda_h + |\nabla^2 w_{h,m,k}|, 0 \}\,dx \leq \gamma + \frac{\delta} {2}.
\end{align}
With \eqref{temp18} at hand, we can repeat the same modification procedure close to $x_1 = a_{h,m,k}$ and obtain,
using a diagonal argument, a sequence $(w_h,\lambda_h) \in \mathcal{A}(\gamma + \delta),
w_h \in W^{2,\infty}(Q,\real^d)$ where $w_h = w_h(\pm a_{h},\cdot)$ near $x_1 = \pm \frac{1}{2}$,
$\nabla w_h = \pm a \otimes e_2$ near $x_2 = \pm \frac{1}{2}$,
$w_h \rightarrow u_0$ in $L^1(Q,\real^d)$ and
$F_{\varepsilon_h} (w_h,\lambda_h) \rightarrow \tilde \Phi(\gamma)$.
We can now reflect $w_h$ with respect to the axis $x_1 = \frac{1}{2}$, translate it such that it is defined
on $(-1,1) \times \big(-\frac{1}{2},\frac{1}{2}\big)$ and denote it by $v_h$. Using property 
(\ref{(H5)}) of $W$, we obtain that we have found a sequence of $2Q'$-periodic functions as desired.
For later use, we also note that $v_h \in W^{2,\infty}_{\text{loc}}
\big(\real \times \big(-\frac{1}{2},\frac{1}{2}\big),\real^d\big)$ if extended periodically.

In order to show (\ref{temp18}), we first claim that
\begin{align} \label{temp19}
	&\limsup_{h \rightarrow \infty} \limsup_{k \rightarrow \infty} \limsup_{m \rightarrow \infty}
		\int_{R_{h,m,k}^-} \frac{1}{\varepsilon_h} W(\nabla w_{h,m,k}) \,dx = 0, \nonumber \\
	&\limsup_{h \rightarrow \infty} \limsup_{k \rightarrow \infty} \limsup_{m \rightarrow \infty}
		\int_{R_{h,m,k}^-} \varepsilon_h | \nabla^2 w_{h,m,k}|^2 \,dx = 0 \quad \text{and} \nonumber \\
	&\limsup_{h \rightarrow \infty} \limsup_{k \rightarrow \infty} \limsup_{m \rightarrow \infty}
		\int_{R_{h,m,k}^-} |w_{h,m,k} - u_h| \,dx = 0.
\end{align}
We are going to exploit that, from the very definition of $w_{h,m,k}$ we have that
\begin{align*}
	\nabla w_{h,m,k} (x) &= \varphi_{h,m,k}(x_1) \nabla u_h(x) + (1 - \varphi_{h,m,k}(x_1))
		\big(  0 \mid \frac{\partial u_h}{\partial x_2}(-a_{h,m,k},x_2) \big) \\
	&+ (u_h(x_1,x_2) - u_h(a_{h,m,k},x_2)) \otimes \varphi_{h,m,k}'(x_1) e_1
\end{align*}
and, by assumption (\ref{(H2)''}), it holds true that
\begin{equation*}
	W(\xi) \leq C(1 + |\xi|^p) \leq C(1 + 2^{p-1}|\nabla u_0|^p + 2^{p-1}|\xi - \nabla u_0|^p ).
\end{equation*}
As a consequence of that, using (\ref{temp15}), (\ref{temp16}) and (\ref{temp17}), we get
\begin{align*}
	\int_{R_{h,m,k}^-}& \frac{1}{\varepsilon_h} W(\nabla w_{h,m,k})\,dx \leq
		\frac{C}{\varepsilon_h} \int_{R_{h,m,k}^-} 1 + |\nabla u_0|^p 
		+ |\nabla w_{h,m,k} - \nabla u_0|^p \,dx \\
	&\leq \frac{C}{\varepsilon_h mk} (1+ |a \otimes e_2|^p)
		+ \frac{C}{\varepsilon_h} \int_{R_{h,m,k}^-} |\nabla u_h - \nabla u_0|^p 
		+ \left| \nabla u_0 - \left(  0 \mid \frac{\partial u_h}{\partial x_2}(-a_{h,m,k},x_2) \right) 
		\right|^p \,dx \\
	&\quad + \frac{C}{\varepsilon_h} \int_{R_{h,m,k}^-} m^p k^p |u_h(x) - u_h(-a_{h,m,k},x_2)|^p \,dx \\
	&\leq \frac{C}{\varepsilon_h mk} 
		+ \frac{C}{\varepsilon_h} \int_{R_{h,m,k}^-} |\nabla u_h - \nabla u_0|^p \,dx
		+ \frac{C}{\varepsilon_h mk} \int_{-\frac{1}{2}}^{\frac{1}{2}}
		|\nabla (u_0 - u_h)(-a_{h,m,k}, x_2)|^p \,dx_2 \\
	& \quad+ \frac{C m^p k^p}{\varepsilon_h}
		\int_{R_{h,m,k}^-} |u_h(x) - u_h(-a_{h,m,k},x_2)|^p \,dx \\
	&\leq \frac{C}{\varepsilon_h mk} + \frac{3 C \tilde \Phi(\gamma)}{\varepsilon_h m^{p+1} k}
		+ \frac{6 C \tilde \Phi(\gamma)}{\varepsilon_h m^{p+1} k}
		+ \frac{C m^p k^p}{\varepsilon_h}
		\int_{R_{h,m,k}^-} |u_h(x) - u_h(-a_{h,m,k},x_2)|^p \,dx.
\end{align*}
By H\"older's inequality we infer that
\begin{align*}
	|u_h(x) - u_h(-a_{h,m,k},x_2)|^p &\leq \left( \int_{-b_{h,m,k}}^{-c_{h,m,k}}
		\left| \frac{\partial u_h}{\partial x_1}(s,x_2) \right| ds \right)^p \\
	&\leq \frac{C}{(mk)^{p-1}}  \int_{-b_{h,m,k}}^{-c_{h,m,k}}
		\left| \frac{\partial u_h}{\partial x_1}(s,x_2) \right|^p ds\\&
		\leq \frac{C}{(mk)^{p-1}}  \int_{-b_{h,m,k}}^{-c_{h,m,k}}
		|\nabla (u_h - u_0)(s,x_2)|^p ds.
\end{align*}
The latter estimate together with \eqref{temp16} gives
\begin{align} \label{temp20}
	\int_{R_{h,m,k}^-} |u_h(x) - u_h(-a_{h,m,k},x_2)|^p \,dx &\leq  
		\frac{C}{(mk)^{p-1}} \int_{R_{h,m,k}^-} \int_{-b_{h,m,k}}^{-c_{h,m,k}}
		|\nabla (u_h - u_0)(s,x_2)|^p ds\,dx \nonumber \\
	&\leq \frac{C}{(mk)^p} \int_{R_{h,m,k}^-} |\nabla u_h - \nabla u_0|^p \,dx
		\leq \frac{C}{m^{2p+1} k^{p+1}}
\end{align}
which shows the first equation in \eqref{temp19}.
To show the second equation, we remark that, by the very definition of $\varphi_{h,m,k}$, in $R_{h,m,k}^-$ we have that
\begin{align*}
	|\nabla^2 w_{h,m,k}(x)|^2 \leq\, &C \left| \nabla^2 u_h \right| ^2 +
		C\big|\frac{\partial^2 u_h}{\partial x_2^2}(-a_{h,m,k},x_2)\big|^2 
	+ Cm^4 k^4 |u_h(x) - u_h(-a_{h,m,k},x_2)|^2 \\&+ Cm^2 k^2 
		\big| \nabla u_h(x) - \left( 0 \mid \frac{\partial u_h}{\partial x_2}(-a_{h,m,k},x_2) \right)
		\big|^2.
\end{align*}
Hence, by \eqref{temp16} and \eqref{temp17} it follows that
\begin{align*}
	\int_{R_{h,m,k}^-} \varepsilon_h |\nabla^2 w_{h,m,k}|^2 \,dx
		&\leq C \int_{R_{h,m,k}^-} \varepsilon_h |\nabla^2 u_h|^2 \,dx
		+ \frac{C \varepsilon_h}{mk} \int_{-\frac{1}{2}}^{\frac{1}{2}}
		|\nabla^2 u_h(-a_{h,m,k},x_2)|^2 \,dx_2 \\
	&\quad + C \varepsilon_h m^4 k^4 \int_{R_{h,m,k}^-} |u_h(x) - u_h(-a_{h,m,k},x_2)|^2 \,dx \\
	&\quad + C \varepsilon_h m^2 k^2 \int_{R_{h,m,k}^-}
		\left| \nabla u_h(x) - \left( 0 \mid \frac{\partial u_h}{\partial x_2}(-a_{h,m,k},x_2) \right)
		 \right|^2 \,dx \\
	&\leq \frac{C}{mk} + C \varepsilon_h m^4 k^4 \int_{R_{h,m,k}^-} |u_h(x) - u_h(-a_{h,m,k},x_2)|^2 
		\,dx \\
	&\quad+ C \varepsilon_h m^2 k^2 \int_{R_{h,m,k}^-}
		\left| \nabla u_h(x) - \left( 0 \mid \frac{\partial u_h}{\partial x_2}(-a_{h,m,k},x_2) \right)
		 \right|^2 \,dx.
\end{align*}
We are left with an estimate of the two integral terms in the previous expression. The first can be estimated noting that H\"older's inequality and \eqref{temp20} yield
\begin{align} \label{temp21}
	\int_{R_{h,m,k}^-} \hspace{-.5cm}|u_h(x) - u_h(-a_{h,m,k},x_2)|^2 \,dx &\leq |{R_{h,m,k}^-}|^{(p-2)/p}
		\left( \int_{R_{h,m,k}^-} \hspace{-.5cm}|u_h(x) - u_h(-a_{h,m,k},x_2)|^p \,dx \right)^{2/p}\nonumber  \\
	&\leq \frac{1}{(mk)^{(p-2)/p}} \left( \frac{C}{m^{2p+1} k^{p+1}} \right)^{2/p}
		= \frac{C}{m^5 k^3}.
\end{align}
In order to estimate the second term we first observe that the following inequality holds true
\begin{align*}
	\Big| \nabla u_h(x) - \Big( 0 \mid \frac{\partial u_h}{\partial x_2}(-a_{h,m,k},x_2) \Big) \Big|^2
		&= \left| \frac{\partial u_h}{\partial x_1}(x) \right|^2 + 
		\left| \frac{\partial u_h}{\partial x_2}(x) - \frac{\partial u_h}{\partial x_2}(-a_{h,m,k},x_2) \right|^2 \\
	&\leq |\nabla u_h(x) - \nabla u_0(x)|^2 + \Big( \int_{-c_{h,m,k}}^{-b_{h,m,k}} 
		\Big|\frac{\partial^2 u_h}{\partial x_1x_2}(s,x_2)\Big|ds \Big)^2 \\
	&\leq |\nabla u_h(x) - \nabla u_0(x)|^2 + \frac{1}{mk}
		\int_{-c_{h,m,k}}^{-b_{h,m,k}} |\nabla^2 u_h(s,x_2)|^2 ds.
\end{align*}
Combining it with \eqref{temp16}, we find that 
\begin{align*}
	\int_{R_{h,m,k}^-} &\left| \nabla u_h(x) - 
		\left( 0 \mid \frac{\partial u_h}{\partial x_2}(-a_{h,m,k},x_2) \right) \right|^2 \,dx \\
	&\leq \int_{R_{h,m,k}^-} |\nabla u_h - \nabla u_0|^2 \,dx + \frac{1}{mk}
		\int_{R_{h,m,k}^-} \int_{-c_{h,m,k}}^{-b_{h,m,k}} |\nabla^2 u_h(s,x_2)|^2 ds\,dx \\
	&\leq \frac{1}{(mk)^{(p-2)/p}} \left( \int_{R_{h,m,k}^-} |\nabla u_h - \nabla u_0|^p \,dx \right)^{2/p}
		+ \frac{1}{m^2 k^2} \int_{R_{h,m,k}^-} |\nabla^2 u_h|^2 \,dx \\
	&\leq \frac{C}{m^3 k} + \frac{C}{\varepsilon_h m^3 k^3}.
\end{align*}
This completes the proof of the second equation in \eqref{temp19},
The third equation follows noting that from \eqref{temp21} it holds that
\begin{equation*}
	\int_{R_{h,m,k}^-} |w_{h,m,k} - u_h|\,dx \leq \frac{1}{mk}
	\left( \int_{R_{h,m,k}^-} |u_h(x) - u_h(-a_{h,m,k},x_2)|^2 \,dx \right)^{1/2}
	\leq \frac{C}{m^{7/2} k^{5/2}}.
\end{equation*}
We now prove \eqref{temp18}. We first observe that from Remark \ref{R6.1} it follows that
\begin{align*}
	W \Big( 0,\frac{\partial u_h}{\partial x_2}(-a_{h,m,k},x_2) \Big) &\leq C\,W(\nabla u_h(-a_{h,m,k},x_2))
		+ C \left| \frac{\partial u_h}{\partial x_1}(-a_{h,m,k},x_2) \right| ^p \\
	&\leq C\,W(\nabla u_h(-a_{h,m,k},x_2)).
\end{align*}
Set $J := (-\frac{1}{2}, -b_{h,m,k}) 
\times (-\frac{1}{2}, \frac{1}{2})$, the previous estimate together with \eqref{temp17}, gives
\begin{align}\label{est:step1}
	F_{\varepsilon_h}&(w_{h,m,k},\lambda_h,J \backslash R_{h,m,k}) \leq
		\int_{J \backslash R_{h,m,k}} \frac{1}{\varepsilon_h}W(\nabla w_{h,m,k}) 
		+ 2\varepsilon_h |\nabla^2 w_{h,m,k}|^2 \,dx \\
	\nonumber&\leq \frac{2}{m} \int_{-\frac{1}{2}}^{\frac{1}{2}} \frac{1}{\varepsilon_h}
		W \left( 0,\frac{\partial u_h}{\partial x_2}(-a_{h,m,k},x_2) \right) 
		+ \varepsilon_h \left| \frac{\partial^2 u_h}{\partial x_2^2}(-a_{h,m,k},x_2) \right|^2 \,dx \\
	\nonumber&\leq \frac{C}{m} \int_{-\frac{1}{2}}^{\frac{1}{2}} \frac{1}{\varepsilon_h}
		W(\nabla u_h(-a_{h,m,k},x_2)) 
		+ \varepsilon_h \left| \frac{\partial^2 u_h}{\partial x_2^2}(-a_{h,m,k},x_2) \right|^2 \,dx
		\leq \frac{C}{m},
\end{align}
This shows that
\begin{equation*}
	F_{\varepsilon_h}(w_{h,m,k},\lambda_h) \leq F_{\varepsilon_h}(u_h,\lambda_h, Q \backslash J)
	+F_{\varepsilon_h}(w_{h,m,k},\lambda_h,R_{h,m,k}^-) + \frac{C}{m}
\end{equation*}
which, combined with \eqref{temp19}, implies the first equation in \eqref{temp18}.
Again by \eqref{temp17} it holds that
\begin{align*}
	\int_Q |w_{h,m,k} - u_0| \,dx \leq& \int_{Q \backslash J} |u_h - u_0|\,dx 
		+ \int_{R_{h,m,k}^-} |w_{h,m,k} - u_0| \,dx \\
	&+ \frac{1}{m} \int_{-\frac{1}{2}}^{\frac{1}{2}} |u_h(-a_{h,m,k},x_2) - u_0(x)| \,dx \\
	\leq &\int_{Q \backslash J} |u_h - u_0|\,dx + \int_{R_{h,m,k}^-} |w_{h,m,k} - u_0| \,dx
		+ \frac{C}{m},
\end{align*}
which, together with \eqref{temp19}, implies the second equation in \eqref{temp18}.
To show the third equation, we first note that by the previous inequalities leading to the proof of \eqref{temp19}, we have that
\begin{equation*}
	\int_{R_{h,m.k}^-} |\nabla^2 w_{h,m,k}|^2 \,dx \leq \frac{C}{\varepsilon_h mk} + \frac{Ck}{m}.
\end{equation*}
From that and \eqref{est:step1} it we get 
\begin{align*}
	\int_J &\max \{ \lambda_h + |\nabla^2 w_{h,m,k}|, 0 \} \,dx \leq \int_J |\nabla^2 w_{h,m,k}| \,dx \\
	&\leq \frac{1}{m} \int_{-\frac{1}{2}}^{\frac{1}{2}} 
		\left| \frac{\partial^2 u_h}{\partial x_2^2}(-a_{h,m,k},x_2) \right|
		+ |R_{h,m,k}^-|^{1/2} \left( \int_{R_{h,m.k}^-} |\nabla^2 w_{h,m,k}|^2 \,dx \right) ^{1/2} \\
	&\leq \frac{C}{\varepsilon_h m} + \left( 
		\frac{C}{\varepsilon_h m^2 k^2} + \frac{C}{m^2} \right)^{1/2}.
\end{align*}
This completes the first part of the proof.
\item[Step 2.]
For $N \geq 3$ one can repeat the argument in step 1, leading to the definition of the sequence of functions $(w_h)$ (from which the sequence $(u_h)$ is obtained), in each coordinate direction. We give here only the main idea (see step 2 in the proof of Proposition 6.3 \cite{cfl} for additional details). Starting from a sequence $(u_h) \subset W^{2,2}(Q,\real^d)$, we modify $u_h$ as in step 1
and obtain, by reflection with respect to the hyperplane $\{ x_1 = \frac{1}{2} \}$, the functions $w_h \in W^{2,\infty}((-1,1) \times (-\frac{1}{2},\frac{1}{2})^{N-1}, \real^d)$ that are $(-1,1)$-periodic in $x_1$. The desired sequence $(w_h)$ is then obtained repeating the same construction with respect to the variables $x_2, \dots, x_{N-1}$.
\end{description}\end{proof}

We define for $\gamma \geq 0$
\begin{align}\label{Phiptilde}
	\tilde \Phi_p(\gamma) := \inf \big\{ &F_{1/L}(u,\lambda) : 
		L>0, u \in W^{2,\infty}(Q,\real^d), \nabla u = \pm a \otimes e_N
		\text{ near } x_N = \pm \frac{1}{2}, \\
	&u \text{ periodic of period one in }x',
		(u,\lambda) \in \mathcal{A}(\gamma) \big\}
\end{align}
and
\begin{equation}\label{Phip}
	\Phi_p(\gamma) := \lim_{\delta \rightarrow 0^+} 
	\tilde \Phi_p(\gamma + \delta).
\end{equation}

\begin{prop}\label{Phi-Phip}
Let  $W:\real^{d\times N} \rightarrow [0,\infty)$ satisify (\ref{(H1)}), (\ref{(H2)''}) , (\ref{(H4)}) and (\ref{(H5)}). Then
it holds that $\Phi_p(\gamma) = \Phi(\gamma)$.
\end{prop}
\begin{proof}
From the propositions above it follows that $\Phi_p \leq \Phi$.
The other inequality can be shown by an application of the same rescaling argument
we are going to use in the proof of Theorem \ref{recoverySequence}, which allows us to define from an admissible pair
$(u,\lambda)$ for $\Phi_p$ an admissible sequence for $\Phi$. 
\end{proof}

In analogy with \cite[Proposition 5.3]{cfl}
If we replace the assumption on the potential (\ref{(H5)}) with the following
\begin{equation*} \tag{H5} \label{(H1d)}
	W(\xi) \geq W(0,\xi_N) , \quad \xi = (\xi',\xi_N) \in \real^{d \times N}
\end{equation*}
we can show that the sequence of functions from \cref{verticalMatching} can be chosen to 
depend only on the $x_N$-variable. More precisely, it holds
\begin{prop} \label{1dProfile}
Let $W:\real^{d\times N} \rightarrow [0,\infty)$ satisfy \eqref{(H1)}, \eqref{(H2)''} , \eqref{(H4)} and
\eqref{(H1d)} and let
$\gamma \geq 0$ and $\delta > 0$ be given. Then there exist sequences
$\varepsilon_h^\delta \rightarrow 0^+$ and $(u_h^\delta, \lambda_h^\delta)
\in \mathcal{A} (\gamma + \delta)$ satisfying
\begin{align*}
	&u_h^\delta \in W^{2,\infty}(Q,\real^d), \\
	& u_h^\delta \rightarrow |x_N|a \text{ in } L^1 (Q,\real^d), \\
	& \nabla u_h^\delta = \pm a \otimes e_N \text{ near } x_N 
		= \pm \frac{1}{2}, \\
	& u_h \text{ depends only on } x_N, \\
	&\lim_{h \rightarrow \infty} F_{\varepsilon_h}
		(u_h^\delta, \lambda_h^\delta) = \tilde \Phi(\gamma).
\end{align*}
\end{prop}
\begin{proof}
As in the proof of \cref{transversalPeriodicity} we obtain the existence of
sequences $\varepsilon_h \rightarrow 0^+, (u_h, \lambda_h) \in \mathcal{A}(\gamma + \delta)$
and $u_h \in C^2(Q,\real^d)$
such that $u_h \rightarrow u_0$ in $L^1(Q,\real^d)$,
$\nabla u_h = \pm a \otimes e_N$ near $x_2 = \pm \frac{1}{2}$ and
$\lim_{h \rightarrow \infty} F_{\varepsilon_h} (u_h,\lambda_h) = \tilde \Phi(\gamma)$.
Because we have that $(x_N \mapsto u_h(0,x_n), \lambda) \in \mathcal{A}(\gamma + \delta)$,
it is enough to observe that
\begin{align*}
	F_{\varepsilon_h} (u_h,\lambda_h) &= \int_Q\frac{1}{\varepsilon_h} W(\nabla u_h) 
		+ \varepsilon_h |\nabla^2 u_h|^2 + \varepsilon_h \min \{\lambda_h^2, |\nabla^2 u_h|^2 \}\,\,dx \\
	&\geq \int_Q\frac{1}{\varepsilon_h} W(0, \partial_{x_N} u_h) 
		+ \varepsilon_h |\nabla^2 u_h(0,x_N)|^2 + \varepsilon_h \min 
		\{\lambda_h^2, |\nabla^2 u_h(0,x_N)|^2 \}\,\,dx \\
	&= F_{\varepsilon_h} (u_h(0,\cdot),\lambda_h).
\end{align*}
\end{proof}

\section{\(\Gamma\)-convergence}\label{gammalim}

In this section we state and prove the two propositions \ref{gammaLimInf} and \ref{recoverySequence},
which together imply our $\Gamma$-convergence result:

\begin{thm} \label{mainthm}
Let $E_\varepsilon$ be as in (\ref{tempKorr4}), where $W: \real^{N \times d} \rightarrow [0,\infty)$
is a continuous double-well potential as in (\ref{tempKorr5}) that satisfies growth conditions as in (\ref{(H1)}),
(\ref{(H2)''}) and (\ref{(H4)}) and is even as in (\ref{(H5)}). If the domain $\Omega \subset \real^N$ 
is open, bounded, has a Lipschitz boundary, is simply connected and for all $t \in \real$ it holds that
the section $\{ (x_1, \dots x_N) \in \Omega : x_N = t \}$ is connected, then
in  the space $W^{1,1}(\Omega, \real^d) \times \mathcal{M}(\Omega)$ it holds that
\begin{equation*}
	\gammalim_{\varepsilon \rightarrow 0^+} E_\varepsilon(u,\mu) = E(u,\mu).
\end{equation*}
Here, we have written
\begin{align}\label{limenergy}
	E(u,\mu) := 
	\begin{cases} 
		\int_{S_{\nabla u}} \Phi \left( \frac{d \mu}{d \HN-1 \llcorner S_{\nabla u}} \right) ,&
		u\in BV(\Omega, \{ A,B \}),\\
	+\infty, &\text{otherwise in }W^{1,1}(\Omega;\real^d)\times\mathcal M(\Omega),
	\end{cases}
\end{align}
with $\Phi$ a nonnegative nonincreasing function as in (\ref{tempKorr6}).
\end{thm}
\begin{proof}
The proof follows from \cref{gammaLimInf} and \cref{recoverySequence}.
\end{proof}

\begin{thm} \label{mainthm-1d}
Let $E_\varepsilon$ be as in (\ref{tempKorr4}), where $W: \real^{N \times d} \rightarrow [0,\infty)$
is a continuous double-well potential as in (\ref{tempKorr5}) that satisfies growth conditions as in (\ref{(H1)}),
(\ref{(H2)''}) and (\ref{(H4)}) and (\ref{(H1d)}). If the domain $\Omega \subset \real^N$ 
is open, bounded, has a Lipschitz boundary, is simply connected and for all $t \in \real$ it holds that
the section $\{ (x_1, \dots x_N) \in \Omega : x_N = t \}$ is connected, then
in  the space $W^{1,1}(\Omega, \real^d) \times \mathcal{M}(\Omega)$ it holds that
\begin{equation*}
	\gammalim_{\varepsilon \rightarrow 0^+} E_\varepsilon(u,\mu) = E(u,\mu),
\end{equation*}
where $E(u,\mu)$ is given in \eqref{limenergy}. The surface tension $\Phi$ in Definition \ref{definitionPhi} is obtained further restricting the admissible set of functions in the cell-problem formula to one-dimensional profiles $u_h(x)=u_h(x_N)$. 
\end{thm}
\begin{proof}
The proof of the statement follows from the proof of Theorem \ref{mainthm} taking into account Proposition
\ref{1dProfile}.
\end{proof}
\begin{rem}
Note that for potentials that do not satisfy (\ref{(H1d)}) Proposition \ref{1dProfile}, and hence Theorem \ref{mainthm-1d}, are false as shown in
\cite[Section 8]{cfl}.
\end{rem}

\begin{prop}[$\Gamma$-liminf inequality]\label{gammaLimInf}
Let $W: \real^{d \times N} \rightarrow [0,\infty)$ 
satisfies (\ref{(H1)}) and (\ref{(H2)''}), $\varepsilon_h \rightarrow 0^+$, and let 
$(u_h)$ be a sequence in $W^{2,2}(\Omega, \real^d)$ and $\mu_h = \rho_h \,dx$ be a
sequence in $\mathcal{M}(\Omega)$ such that 
\begin{equation*}
\varepsilon_h \rightarrow 0^+,  \quad u_h \rightarrow u \text{ in }
L^1(\Omega, \real^d), \quad \mu_h \weak* \mu \in \mathcal{M}(\Omega).
\end{equation*}
Then it holds that
\begin{equation*}
E(u,\mu) \leq \liminf_{h \rightarrow \infty} E_{\varepsilon_h}(u_h,\rho_h).
\end{equation*}
\end{prop}
\begin{proof}
Thanks to \eqref{compactness}, up to choosing a subsequence, we may assume that
\begin{align*} &\liminf_{h \rightarrow \infty} E_{\varepsilon_h}(u_h,\rho_h)
	= \lim_{h \rightarrow \infty} E_{\varepsilon_h}(u_h,\rho_h), \\
&\nabla u_h \rightarrow \nabla u \text{ a.e. and } 
	\nabla u \in \BV(\Omega, \{A,B\})  \quad \text{and} \\
&\big ( \frac{1}{\varepsilon_h}W(\nabla u_h) + \varepsilon_h(|\nabla^2u_h|^2 
	+(\rho_h-|\nabla^2u_h|)^2 ) \big ) \mathcal{L}^N \weak* \sigma.
\end{align*}
To conclude it is enough to show the following claim:
\begin{equation*}
\frac{d \sigma}{d {\mathcal H}^{N-1} \llcorner S_{\nabla u}} \geq \Phi 
\Big( \frac{d \mu}{d {\mathcal H}^{N-1} \llcorner S_{\nabla u}} \Big) \quad {\mathcal H}^{N-1}
\text{-a.e. in } S_{\nabla u}.
\end{equation*}
We can use \eqref{laminates} together with 
the Besikovitch differentiation theorem (\cite[Theorem 2.22]{afp})
and \cite[Proposition 1.62]{afp} to see that for ${\mathcal H}^{N-1}$-a.e. 
$x_0 \in S_{\nabla u}$ and for all but at most countably many
$\delta \ll 1$ it holds that
\begin{align*}
&\mu(\partial Q_{x_0,\delta}) = 0 \Rightarrow \lim_{h \rightarrow \infty}
	\int_{Q_{x_0,\delta}} \rho_h \,dx = \mu(Q_{x_{0,\delta}}) \\
&\sigma(\partial Q_{x_0,\delta}) = 0 \Rightarrow \lim_{h \rightarrow \infty}
	E_{\varepsilon_h}(u_h,\rho_h,Q_{x_0,\delta}) = \sigma(Q_{x_{0,\delta}}) \\
&{\mathcal H}^{N-1} \llcorner S_{\nabla u}(Q_{x_0,\delta}) = \delta^{N-1} \\
&\lim_{\delta \rightarrow 0^+} \frac{ \mu(Q_{x_0,\delta})}
	{ {\mathcal H}^{N-1} \llcorner S_{\nabla u} (Q_{x_0,\delta})} =
	\frac{d \mu}{d {\mathcal H}^{N-1} \llcorner S_{\nabla u}}(x_0) =: \rho \\
&\lim_{\delta \rightarrow 0^+} \frac{ \sigma(Q_{x_0,\delta})}
	{ {\mathcal H}^{N-1} \llcorner S_{\nabla u} (Q_{x_0,\delta})} =
	\frac{d \sigma}{d {\mathcal H}^{N-1} \llcorner S_{\nabla u}}(x_0).
\end{align*}
In particular, given $r > 0$, we have 
\begin{align} \label{tempKorr2a}
\int_{Q_{x_0,\delta}} \rho_h \,dx &\leq (1+2r)\delta^{N-1} \rho, \\ \label{tempKorr2b}
(1+2r)\frac{d \sigma}{d {\mathcal H}^{N-1} \llcorner S_{\nabla u}}(x_0) &\geq 
	E_{\varepsilon_h}(u_h,\rho_h,Q_{x_0,\delta}) \delta^{1-N}
\end{align}
for all but at most countably many $\delta \ll 1$ and $h$ sufficiently large.
We define $\lambda_h$ such that
\begin{equation} \label{tempKorr3}
\int_{Q_{x_0,\delta}} \max \{ \lambda_h + |\nabla^2 u_h|, 0 \} \,dx
= \int_{Q_{x_0,\delta}} \rho_h \,dx
\end{equation}
and consider the following subsequences (if they exist) of $(\lambda_h)$ that we do not
relabel:
\begin{description}
\item[Case 1 ($\lambda_h$ satisfies $\lambda_h \geq 0$)]
We have by \eqref{tempKorr2a} and \eqref{tempKorr3}
\begin{equation*}
\int_{Q_{x_0,\delta}} |\nabla^2 u_h| \,dx \leq \int_{Q_{x_0,\delta}} \rho_h \,dx
\leq (1+2r)\delta^{N-1} \rho.
\end{equation*}
This means that $(u_h,0) \in \mathcal{A}((1+2r)\rho, Q_{x_0,\delta})$
and we obtain
\begin{equation*}
\lim_{h \rightarrow \infty} E_{\varepsilon_h}(u_h,\rho_h,Q_{x_0,\delta})
\geq \liminf_{h \rightarrow \infty} F_{\varepsilon_h}(u_h,0,Q_{x_0,\delta})
\geq \tilde \Phi((1+2r)\rho) \delta^{N-1}
\end{equation*}
where in the last inequality we have used \cref{L4.3}.
By \eqref{tempKorr2b} we obtain that
\begin{equation*}
(1+2r)\frac{d \sigma}{d {\mathcal H}^{N-1} \llcorner S_{\nabla u}}(x_0) \geq \tilde \Phi((1+2r)\rho)
\end{equation*}
and, letting $r \rightarrow 0^+$,
\begin{equation*}
\frac{d \sigma}{d {\mathcal H}^{N-1} \llcorner S_{\nabla u}}(x_0) \geq \Phi(\rho).
\end{equation*}
\item[Case 2 ($\lambda_h$ satisfies $\lambda_h \leq 0$)]
By (\ref{tempKorr3}) we have $(u_h,\lambda_h) \in \mathcal{A}((1+2r)\rho, Q_{x_0,\delta})$ and,
since $\lambda_h \leq 0$,
\begin{equation*}
\int_{Q_{x_0,\delta}} \rho_h \,dx \leq \int_{Q_{x_0,\delta}} |\nabla^2 u_h| \,dx
\end{equation*}
and we can apply Lemma \ref{L3.2} to get
\begin{align*}
E_{\varepsilon_h}(u_h,\rho_h,Q_{x_0,\delta}) &=
	\int_{Q_{x_0,\delta}} \frac{1}{\varepsilon_h} W(\nabla u_h) + \varepsilon_h
	|\nabla^2 u_h|^2 + \varepsilon_h (\rho_h - |\nabla^2 u_h|)^2 \,dx \\
&\geq \int_{Q_{x_0,\delta}} \frac{1}{\varepsilon_h} W(\nabla u_h) 
	+ \varepsilon_h |\nabla^2 u_h|^2 
	+ \varepsilon_h \min \{ \lambda_h^2, |\nabla^2 u_h|^2 \} \,dx \\
&= F_{\varepsilon_h}(u_h,\lambda_h,Q_{x_0,\delta}).
\end{align*}
This implies
\begin{equation*}
\lim_{h \rightarrow \infty} E_{\varepsilon_h}(u_h,\rho_h,Q_{x_0,\delta})
\geq \tilde \Phi((1+2r) \rho) \delta^{N-1}.
\end{equation*}
\end{description}
We conclude as in the previous case.
\end{proof}

From now on, we assume that  $\Omega$ is open, bounded, with Lipschitz boundary,
and, in addition, we assume it to be simply connected.
We write $\Omega_t := \{(x',x_N) \in \Omega\ | x_N=t \}$ and set
\begin{equation*}
	\alpha := \inf\{ x_N : \Omega_{x_N} \neq \emptyset \}, \quad
	\beta := \sup\{ x_N : \Omega_{x_N} \neq \emptyset \}.
\end{equation*}
We assume moreover that the sets $\Omega_{x_N}$ are connected for any
$x_N \in (\alpha,\beta)$.\\

In what follows, given a set $A\subset\real^n$ we define $A_\delta := \{ x \in \real^n | d(x,A) < \delta \}$. 
For the sake of simplicity we also use the notation $\Omega_{x_N,\delta}:=(\Omega_{x_N})_\delta$.

\begin{prop}[$\Gamma$-limsup inequality] \label{recoverySequence}
Let $W: \real^{d \times N} \rightarrow [0,\infty)$ satisfy (\ref{(H1)}), (\ref{(H2)''}),
(\ref{(H4)}) and (\ref{(H5)}).
Given a measure $\mu \in \mathcal{M}(\Omega)$
and $u \in W^{1,1}(\Omega,\real^d)$ with $\nabla u \in \BV(\Omega,\{A,B\})$,
and given any sequence $\varepsilon_h \rightarrow 0^+$,
there exist sequences $(\mu_h) \subset \mathcal{M}(\Omega)$
and $(u_h) \subset W^{2,2}(\Omega,\real^d)$ satisfying
\[ \mu_h \weak* \mu, \quad u_h \rightarrow u \text{ in } W^{1,1}(\Omega,\real^d),
\quad \limsup_{h \rightarrow \infty} E_{\varepsilon_h}(u_h,\mu_h)
\leq E(u,\mu). \]
\end{prop}
\begin{proof}
We will prove the statement in several steps. At each step we assume $u$ and $\mu$ to be of increasing generality and provide for them a recovery sequence.
\begin{description}
\item[Step 1.1]
	We assume that $u$  and $\mu$ are given by
	\begin{equation} \label{tempKorr7}
		u(x)=|x_N|a \quad \text{and} \quad  
			\mu = \gamma_1 \chi_K {\mathcal H}^{N-1} \llcorner S_{\nabla u} + \beta \delta_{x_0},
	\end{equation}
	where $\gamma_1, \beta \geq 0$, $x_0 \in \Omega \backslash S_{\nabla u}$ and
	$K \subset \Omega$
	is a compact set. Note that in this case there is only one connected interface, namely $S_{\nabla u}=\Omega_0$.
	We choose $h > 0$ such that $[-4h,4h] \subset [\alpha,\beta]$
	and assume moreover without loss of generality that
	$K = K' \times \big[ -\frac{h}{2},\frac{h}{2} \big]$, with $K' \subset \Omega_0$ compact.
	Given $\varepsilon,\delta, \tilde \delta > 0$, thanks to Proposition \ref{transversalPeriodicity}, the very definition of $\tilde \Phi_p$ and $\Phi_p$ and Proposition \ref{Phi-Phip}, we can choose $L_1>0$ and
	\begin{align*}
		&(v^1,\lambda^1) \in \mathcal{A}(\gamma_1 + \tilde \delta) 
			\cap W^{2,\infty}(Q,\real^d) \\
		&\nabla v^1 = \pm a \otimes e_N
			\text{ near } x_N = \pm \frac{1}{2},\,\,
			v^1 \text{ periodic of period one in } x'
	\end{align*}
	that we assume to be extended periodically in $x'$ such that
	\begin{equation*}
	F_{1/L^1}(v^1,\lambda^1) \leq  \tilde \Phi_p(\gamma_1 + \tilde \delta) 
	+ \tilde \delta \leq \Phi(\gamma_1) + \tilde \delta.
	\end{equation*}
	We set
	\begin{equation*}
		z^1_\varepsilon(x) := \left\{
		\begin{array}{rll}
		&\varepsilon L^1 v^1(\frac{x'}{\varepsilon L^1},-\frac{1}{2})
			- a(x_N + \frac{\varepsilon L^1}{2}), &
			x_N < -\frac{\varepsilon L^1}{2} \\
		&\varepsilon L^1 v^1(\frac{x}{\varepsilon L^1}), &
			|x_N| \leq \frac{\varepsilon L^1}{2} \\
		&\varepsilon L^1 v^1(\frac{x'}{\varepsilon L^1},\frac{1}{2})
			+ a(x_N - \frac{\varepsilon L^1}{2}), &
			x_N > \frac{\varepsilon L^1}{2} \\
		\end{array} \right.
	\end{equation*}
	and
	\begin{equation*}
		\lambda^1_\varepsilon := \frac{\lambda^1}{\varepsilon L^1}.
	\end{equation*}
	 It holds that
	 \begin{equation*}
		\nabla z^1_\varepsilon(x) = \left\{
		\begin{array}{rll}
		&-a \otimes e_N, & x_N < - \frac{\varepsilon L^1}{2} \\
		&\nabla v^1(\frac{x}{\varepsilon L^1}), & |x_N|
			 < \frac{\varepsilon L^1}{2} \\
	 	&a \otimes e_N, & x_N > \frac{\varepsilon L^1}{2}. \\
	 	\end{array} \right.
	\end{equation*}
	Moreover, by the Riemann-Lebesgue lemma, we also have that
	\begin{align*}
		\int_Q \max\{ |\nabla^2 z^1_\varepsilon| +
			\lambda^1_\varepsilon, 0 \} \,dx
			&= \int_{Q'} \int_{-1/2}^{1/2} \max\Big\{ \Big|\nabla^2 v^1
			\Big(\frac{x'}{\varepsilon L^1},x_N\Big)\Big| + \lambda^1, 0 \Big\} \,dx_N \,dx' \\
		&\rightarrow_{\varepsilon \rightarrow 0^+} \int_Q
			\max\{ |\nabla^2 v^1| + \lambda^1, 0 \} \,dx \leq \gamma_1 
			+ \tilde\delta,
	\end{align*}
	where the last inequality follows by the assumption on $(v^1,\lambda^1)$.
	We repeat the same construction, up to replacing the index $1$ by $2$, for $\gamma_2 := 0$ 
	and obtain $(z^2_\varepsilon, \lambda^2_\varepsilon)$.
	We choose a smooth cut-off function $\varphi_\delta: \Omega_{0,2\delta}
	\times (-\frac{h}{2},\frac{h}{2}) \rightarrow \real$ satisfying
	\begin{equation*}
		\varphi_\delta |_K = 1, \quad \varphi_\delta |_{\Omega_{0,2\delta}
		\times (-\frac{h}{2},\frac{h}{2}) \backslash K_\delta} = 0,
		\quad |\nabla \varphi_\delta| \leq \frac{C}{\delta}, \quad
		|\nabla^2 \varphi_\delta| \leq \frac{C}{\delta^2}
	\end{equation*}
	and set $z_{\varepsilon,\delta} := \varphi_\delta z_\varepsilon^1 +
	(1-\varphi_\delta) z_\varepsilon^2$.
	We also choose a smooth cut-off function $\psi_\delta:\Omega \rightarrow \real$
	satisfying
	\begin{equation*}
		\psi_\delta |_{\Omega_{0,\delta} \times (-h/3,h/3)} = 1, 
		\quad \psi_\delta |_{\Omega \backslash \Omega_{0,2\delta} 
		\times (-h/2,h/2)}= 0,
		\quad |\nabla \psi_\delta| \leq \frac{C}{\delta}, \quad
		|\nabla^2 \psi_\delta| \leq \frac{C}{\delta^2}
	\end{equation*}
	and set
	\begin{equation*}
		u_{\varepsilon,\delta} := \psi_\delta z_{\varepsilon,\delta}
		+ (1 - \psi_\delta) u.
	\end{equation*}
	An example of this situation is shown  in Fig.\ref{fig1}.
	\begin{figure}
\begin{tikzpicture}[scale=0.85]
	\filldraw[fill=black!5] (0.5,-3) rectangle (7,3);
	\filldraw[fill=black!12] (1,-2) rectangle (6.5,2);
	\filldraw[fill=black!19] (2,-2) rectangle (3.3,2);
	\draw (3,1) node[anchor=north east]{$K$};
	\filldraw[fill=black!12] (2,2) rectangle (3.3,3);
	\filldraw[fill=black!12] (2,-2) rectangle (3.3,-3);
	
	\draw (7,2) node[anchor=south east]{$\Omega_{0,2\delta}
		\times \left( -\frac{h}{2}, \frac{h}{2} \right)$};
	\draw (6.5,1) 	node[anchor=south east]{$\Omega_{0,\delta}
		\times \left( -\frac{h}{3}, \frac{h}{3} \right)$};
	
	\draw [->](-1,0) -- (9,0) node[anchor=west]{$x_1$};
	\draw [->](0,-5) -- (0,5) node[anchor=south]{$x_2$};
	\draw (-2pt,2) -- (2pt,2) node[anchor=east]{$\frac{h}{3}$};
	\draw (-2pt,3) -- (2pt,3) node[anchor=east]{$\frac{h}{2}$};
	\draw (-2pt,-2) -- (2pt,-2) node[anchor=east]{$-\frac{h}{3}$};
	\draw (-2pt,-3) -- (2pt,-3) node[anchor=east]{$-\frac{h}{2}$};
	
	\draw (1.5,3) -- (1.5,-3);
	\draw (1.5,-3) .. controls (1.5,-5) and (3,-4) .. (4,-2);
	\draw (4,-2) parabola bend (6,0) (6,0);
	\draw (7.5,0) parabola (8.5,1);
	\draw (8.5,1) .. controls (9.5,3) and (1.5,5) .. (1.5,3)
		node[anchor = north, near end]{$\Omega$};
\end{tikzpicture}
\caption{Sketch of the construction in Step 1.1}
\label{fig1}
\end{figure}
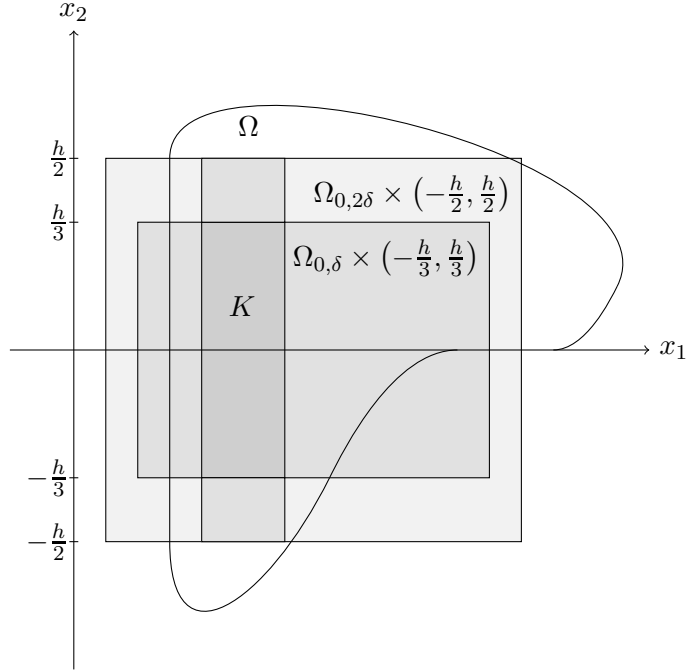 
We claim that, uniformly with respect to $\tilde\delta$,
	\begin{equation*}
	\lim_{\delta \rightarrow 0^+}
	\lim_{\varepsilon \rightarrow 0^+} \Arrowvert u_{\varepsilon,\delta} - u\Arrowvert_{W^{1,p}
	(\Omega,\real^d)} = 0.
	\end{equation*}
	To prove the claim we first observe that
	\begin{equation*}
	|\nabla u_{\varepsilon,\delta} - \nabla u| \leq |\nabla \psi_\delta|
	|z_{\varepsilon,\delta} - u| + |\nabla z_{\varepsilon,\delta} - \nabla u|.
	\end{equation*}
	Since for $|x_N| \geq \max \{ \frac{\varepsilon L^1}{2}, \frac{\varepsilon L^2}{2} \}$
	it holds that $\nabla z_\varepsilon^1 = \nabla z_\varepsilon^2 = \nabla u$, 
	for such $x$ we have that $\nabla z_{\varepsilon,\delta}-\nabla u = 
	\nabla \varphi_\delta (z_\varepsilon^1 - z_\varepsilon^2)$ and,
	recalling that $v^1$ and $v^2$ are bounded in $L^\infty$, it follows that
	\begin{equation*}
		|\nabla u_{\varepsilon,\delta} - \nabla u| \leq |\nabla \psi_\delta||z_{\varepsilon,\delta} - u|
		+ |\nabla z_{\varepsilon,\delta} - \nabla u| =  |\nabla \psi_\delta||z_{\varepsilon,\delta} - u|
		+ |\nabla \varphi_\delta (z_\varepsilon^1 - z_\varepsilon^2)| \leq \frac{C}{\delta} \varepsilon.
	\end{equation*}
	In the set $\{ \psi_\delta = 1 \}$ it holds
	\begin{equation*}
	|\nabla u_{\varepsilon,\delta} - \nabla u| \leq |\nabla \varphi_\delta|
	(|z^1_\varepsilon| + |z^2_\varepsilon|) + |\nabla z^1_\varepsilon|
	+ |\nabla z^2_\varepsilon|
	\leq \frac{C}{\delta}\varepsilon + C
	\end{equation*}
	and on $\{ \psi_\delta = 0 \}$ it holds
	\begin{equation*}
	|\nabla u_{\varepsilon,\delta} - \nabla u| = 0.
	\end{equation*}
	Without loss of generality we may assume that
	$\max \{ \frac{\varepsilon L^1}{2}, \frac{\varepsilon L^2}{2} \} \leq 
	\frac{h}{3}$, which implies that $\varphi_\delta$ is constant on $ \{ |x_N| \leq \max \{ 
	\frac{\varepsilon L^1}{2}, \frac{\varepsilon L^2}{2} \} \}
	\cap \{ 0 < \psi_\delta < 1 \}$ and we have 
	\begin{equation*}
	|\nabla u_{\varepsilon,\delta} - \nabla u| \leq |\nabla \psi_\delta|
	(|z^1_\varepsilon| + |z^2_\varepsilon|) + |\nabla z^1_\varepsilon|
	+ |\nabla z^2_\varepsilon| 
	\leq \frac{C}{\delta}\varepsilon + C.
	\end{equation*}
	It follows
	\begin{equation*}
	\int_\Omega | \nabla u_{\varepsilon,\delta} - \nabla u |^p \,dx
	\leq \frac{C}{\delta^p}\varepsilon^p 
	+ \left| \frac{C}{\delta}\varepsilon + C \right|^p \max \left\{ \frac{\varepsilon L^1}{2},
	 \frac{\varepsilon L^2}{2} \right\} {\mathcal H}^{N-1}(\Omega_{0,2\delta})
	\end{equation*}
	Since $v^{1,2}_{\varepsilon,\delta} \in L^\infty$, it holds that
	$\Arrowvert z_\varepsilon^{1,2} - u \Arrowvert_\infty \leq C \varepsilon$ and the claim follows.
	\newline
	We define a sequence of measures converging to $\mu$.
	We set $\mu^1_\varepsilon := \rho^1_\varepsilon \mathcal{L}^N$, where
	\begin{equation*}
	\rho^1_\varepsilon(x) := \left\{
	\begin{array}{rll}
	&\max\{|\nabla^2 z^1_\varepsilon| + \lambda^1_\varepsilon,0 \}, &
		|x_N| < \frac{\varepsilon L^1}{2} \\
	&\frac{1}{2\sqrt{\varepsilon}} (\gamma_1 + \tilde \delta - \int_{Q}
		\max\{|\nabla^2 v^1| + \lambda^1,0 \}), &
		\frac{\varepsilon L^1}{2} < |x_N|
		< \frac{\varepsilon L^1}{2} + \sqrt{\varepsilon},\\
	&\frac{\beta}{\sqrt{\varepsilon} |B_1|}, &
		x \in B_{\varepsilon^\frac{1}{2N}}(x_0) \\
	&0, & \text{ otherwise}.
	\end{array} \right.
	\end{equation*}
	We have that $\mu^1_\varepsilon \weak* (\gamma_1  +\tilde \delta)
	{\mathcal H}^{N-1} \llcorner S_{\nabla u} + \beta \delta_{x_0}$:
	For any open cylinder $\tilde \Omega \subset\subset \Omega$ of the type $\tilde \Omega
	= \tilde \Omega' \times (b,c)$ with $b < 0 < c$ 
	 we have, by the 
	Riemann-Lebesgue lemma, that
	\begin{align*}
	\int_{\Omega} \chi_{\tilde \Omega} d\mu^1_\varepsilon &=
		\int_{-\varepsilon L^1/2}^{\varepsilon L^1/2} \int_{\tilde \Omega'}
		\max\{|\nabla^2 z^1_\varepsilon| + \lambda^1_\varepsilon,0 \} \,dx \\
	&\quad + \frac{1}{\sqrt{\varepsilon}}
		\int_{\varepsilon L^1/2}^{(\varepsilon L^1/2) + \sqrt{\varepsilon}}
		\int_{\tilde \Omega'} \big( \gamma_1 + \tilde \delta -
		\int_Q \max\{|\nabla^2 v^1| + \lambda^1,0 \}
		\big) \,dx + \beta \chi_{\tilde \Omega}(x_0) \\
	&=\int_{-1/2}^{1/2} \int_{\tilde \Omega'}
		\max \left\{ \left| \nabla^2 v^1\Big(\frac{x'}{\varepsilon L^1},x_N\Big) \right|
		+ \lambda^1,0 \right\} \,dx \\
	&\quad + {\mathcal H}^{N-1}(\tilde \Omega') (\gamma_1 + \tilde \delta
		-\int_Q \max\{|\nabla^2 v^1| + \lambda^1,0 \}) + \beta \chi_{\tilde \Omega}(x_0) \\
	&\rightarrow_{\varepsilon \rightarrow 0^+} {\mathcal H}^{N-1}(\tilde \Omega') (\gamma_1 + \tilde \delta)
		+ \beta \chi_{\tilde \Omega}(x_0)
		= {\mathcal H}^{N-1}(S_{\nabla u} \cap \tilde \Omega) (\gamma_1 + \tilde \delta)
		+ \beta \chi_{\tilde \Omega}(x_0).
	\end{align*}
	If instead $0 \notin (b,c)$ we have that
	\begin{equation*}
	\int_{\Omega} \chi_{\tilde \Omega} d\mu^1_\varepsilon \rightarrow
	\beta \chi_{\tilde \Omega}(x_0).
	\end{equation*}
	To obtain the claimed weak* convergence it is enough to recall that
	compactly supported continuous functions can be approximated 
	uniformly by piecewise constant functions on cylinders of the type of $\tilde \Omega$.
	Finally, one can analogously define $\rho^2_\varepsilon$ for
	$\gamma_2 := 0$.
	Eventually, setting $\rho_{\varepsilon,\delta} := \varphi_\delta \rho^1_\varepsilon
	+ (1- \varphi_\delta) \rho^2_\varepsilon$ and $\mu_{\varepsilon,\delta}:=\rho_{\varepsilon,\delta}\mathcal{L}^N$ we get
	\begin{equation*}
	\mu_{\varepsilon,\delta}
	\weak*_{\varepsilon \rightarrow 0^+} \varphi_\delta (\gamma_1 + \tilde \delta)
	{\mathcal H}^{N-1} \llcorner S_{\nabla u} + \beta \delta_{x_0},
	\end{equation*}
	which implies that $\lim_{\tilde \delta \rightarrow 0^+} \lim_{\delta \rightarrow 0^+}
	\lim_{\varepsilon \rightarrow 0^+} \mu_{\varepsilon,\delta}\mathcal{L^N} 
	= \mu$  with respect to the weak* convergence.
	Let us observe that in particular $\mu_{\varepsilon,\delta}(\Omega)$
	is uniformly bounded in $\varepsilon,  \delta$ and that the bound decreases
	for decreasing $\tilde \delta$.
	Next we claim that
	\begin{equation*}
	\limsup_{\tilde \delta\rightarrow 0^+} \limsup_{\delta \rightarrow 0^+}
	\limsup_{\varepsilon \rightarrow 0^+} E_\varepsilon
	(u_{\varepsilon,\delta}, \rho_{\varepsilon,\delta}) \leq E(u,\mu). 
	\end{equation*}
	We have that
	\begin{equation*}
	|\nabla^2 u_{\varepsilon,\delta}| \leq |\nabla^2 \psi_\delta|
	|z_{\varepsilon,\delta} - u| + 2|\nabla \psi_\delta|
	|\nabla z_{\varepsilon,\delta} - \nabla u| 
	+ |\nabla^2 z_{\varepsilon,\delta}|.
	\end{equation*}
	Since for $|x_N| \geq \max \{ \frac{\varepsilon L^1}{2} \frac{\varepsilon L^2}{2} \}$
	it holds that $\nabla z_\varepsilon^1 = \nabla z_\varepsilon^2 = \nabla u$, 
	for such $x$ and for $\delta$ sufficiently small we obtain
	\begin{equation*}
	|\nabla^2 u_{\varepsilon,\delta}| \leq \frac{C}{\delta^2} \varepsilon.
	\end{equation*}
	In the set $\{ \psi_\delta = 1 \} \cap \{ 0 < \varphi_\delta < 1 \}$,
	recalling that $v^{1,2} \in W^{2,\infty}(\Omega, \real^d)$, for $\delta$ sufficiently small,
	it holds
	\begin{align*}
	|\nabla^2 u_{\varepsilon,\delta}| &= |\nabla^2 z_{\varepsilon,\delta}| 
		\leq |\nabla^2 \varphi_\delta| |z^1_\varepsilon - z^2_\varepsilon| 
		+ 2|\nabla \varphi_\delta||\nabla z^1_\varepsilon - \nabla
		z^2_\varepsilon| + |\nabla^2 z^1_\varepsilon| 
		+ |\nabla^2 z^2_\varepsilon| \\
	&\leq \frac{C}{\delta^2}\varepsilon + \frac{C}{\delta} 
	+ \frac{C}{\varepsilon} \leq \frac{C}{\delta^2} + \frac{C}{\varepsilon}.
	\end{align*}
	On $\{ \psi_\delta = 1 \} \cap \{ \varphi_\delta = 1 \}$ 
	it holds
	\begin{equation*}
	|\nabla u_{\varepsilon,\delta} - \nabla u| = |\nabla z^1_\varepsilon| 
	\leq C, \quad
	\nabla^2 u_{\varepsilon,\delta} = \nabla^2 z^1_\varepsilon
	\end{equation*}
	and similarly on $\{ \psi_\delta = 1 \} \cap \{ \varphi_\delta = 0 \}$. 
	On $\{ \psi_\delta = 0 \}$ it holds $\nabla u_{\varepsilon,\delta} = \nabla u$, hence 
	\begin{equation*}
	|\nabla^2 u_{\varepsilon,\delta}| = 0.
	\end{equation*}
	In the next estimates we may assume that
	$\max \big\{ \frac{\varepsilon L^1}{2}, \frac{\varepsilon L^2}{2} \big\} \leq 
	\frac{h}{3}$. As a result,  $\varphi_\delta$ is constant on $ \big\{ |x_N| \leq \max \big\{ 
	\frac{\varepsilon L^1}{2}, \frac{\varepsilon L^2}{2} \big\} \big\}
	\cap \{ 0 < \psi_\delta < 1 \}$,
	hence
	\begin{align*}
	|\nabla^2 u_{\varepsilon,\delta}| &\leq |\nabla^2 \psi_\delta|
		|z_{\varepsilon,\delta} - u| + 2|\nabla \psi_\delta|
		|\nabla z_{\varepsilon,\delta} - \nabla u| 
		+ |\nabla^2 z_{\varepsilon,\delta}| \\
	&\leq \frac{C}{\delta^2}\varepsilon + \frac{C}{\delta} 
	+ \frac{C}{\varepsilon} \leq \frac{C}{\delta^2} + \frac{C}{\varepsilon}
	\end{align*}
	for $\delta$ sufficiently small.
	As a consequence of the previous estimates, setting $J_\varepsilon:=\{ |x_N| > \max \{ \frac{\varepsilon L^1}{2},
		 \frac{\varepsilon L^2}{2} \} \}$ and $K_\varepsilon:=\{ \frac{\varepsilon L^1}{2} < |x_N| 
		< \frac{\varepsilon L^1}{2} + \sqrt{\varepsilon} \}$ we get
	\begin{align*}
	&\int_{\Omega \cap J_\e}
		 \frac{1}{\varepsilon}|\nabla u_{\varepsilon,\delta} - \nabla u|^p
		 + \varepsilon |\nabla^2 u_{\varepsilon,\delta}|^2
		 + \varepsilon (\rho_{\varepsilon,\delta} - 
		 |\nabla^2 u_{\varepsilon,\delta}|)^2 \,dx \\
	&\quad \leq \frac{C}{\delta^p} \varepsilon^{p-1} 
		+ 3 \frac{C}{\delta^4} \varepsilon^3 	
		+ 2 \varepsilon \int_{\Omega \cap K_\e} 
		\left( \frac{1}{2\sqrt{\varepsilon}} \left(\gamma_1 + \tilde \delta - \int_{Q}
		\max\{|\nabla^2 v^1| + \lambda^1,0 \}\right) \right)^2 \,dx \\
	&\quad \quad+ 2\varepsilon \beta^2
		+ 2 \varepsilon \int_{\Omega \cap K_\e} 
		\left( \frac{1}{2\sqrt{\varepsilon}} \left(\gamma_2 + \tilde \delta - \int_{Q}
		\max\{|\nabla^2 v^2| + \lambda^2,0 \}\right) \right)^2 \,dx \\
	&\quad \leq \frac{C}{\delta^p} \varepsilon^{p-1} 
		+ \frac{C}{\delta^4} \varepsilon^3 +2\varepsilon \beta^2 
		+ C \sqrt{\varepsilon},
	\end{align*}

	\begin{align*}
	&\int_{\{ \psi_\delta = 1 \} \cap \{ 0 < \varphi_\delta < 1 
		\} \cap J_\e}
		\frac{1}{\varepsilon}|\nabla u_{\varepsilon,\delta} - \nabla u|^p
		 + \varepsilon |\nabla^2 u_{\varepsilon,\delta}|^2
		 + \varepsilon (\rho_{\varepsilon,\delta} - 
		 |\nabla^2 u_{\varepsilon,\delta}|)^2 \,dx \\
	&\quad\leq \left|K_\delta \backslash K \cap \left\{ |x_N| 
		<  \max \left\{ \frac{\varepsilon L^1}{2},\frac{\varepsilon L^2}{2}\right\} \right\}\right|
		\left( \frac{1}{\varepsilon}
		\left(\frac{C}{\delta}\varepsilon - C\right)^p
		+ 3 \varepsilon \left(\frac{C}{\delta^2} + \frac{C}{\varepsilon}\right)^2 \right) \\
	&\quad \quad+ 2 \varepsilon \left(\mu_{\varepsilon,\delta} \left(K_\delta \backslash K 
		\cap \left\{ |x_N| <  \max \left\{ \frac{\varepsilon L^1}{2},
		\frac{\varepsilon L^2}{2}\right\} \right\}\right)\right)^2 \\
	&\quad\leq \max \{\varepsilon L^1,\varepsilon L^2 \}
		\delta^{N-1} \left( \frac{1}{\varepsilon}
		\left(\frac{C}{\delta}\varepsilon - C\right)^p
		+ 3 \varepsilon \left(\frac{C}{\delta^2} + \frac{C}{\varepsilon}\right)^2 \right)
		+ \varepsilon C
	\end{align*}
	and
	\begin{equation*}
	\int_{\{ \psi_\delta = 0 \}} 
		\frac{1}{\varepsilon}|\nabla u_{\varepsilon,\delta} - \nabla u|^p
		 + \varepsilon |\nabla^2 u_{\varepsilon,\delta}|^2
		 + \varepsilon (\rho_{\varepsilon,\delta} - 
		 |\nabla^2 u_{\varepsilon,\delta}|)^2 \,dx 
		 \leq \varepsilon C.
	\end{equation*}
	Moreover, recalling that we have assumed that $\max \big\{ \frac{\varepsilon L^1}{2},
	\frac{\varepsilon L^2}{2} \big\} \leq  \frac{h}{3}$, we also get
	\begin{align*}
	&\int_{\{ 0 < \psi_\delta < 1 \} \cap 
		J_\e}
		\frac{1}{\varepsilon}|\nabla u_{\varepsilon,\delta} - \nabla u|^p
		 + \varepsilon |\nabla^2 u_{\varepsilon,\delta}|^2
		 + \varepsilon (\rho_{\varepsilon,\delta} - 
		 |\nabla^2 u_{\varepsilon,\delta}|)^2 \,dx \\
	&\quad\leq \left| (\Omega_{0,2\delta} \backslash \Omega_{0,\delta}) \times 
		\left(-\max \left\{ \frac{\varepsilon L^1}{2},
		\frac{\varepsilon L^2}{2} \right\}, \max \left\{ \frac{\varepsilon L^1}{2},
		\frac{\varepsilon L^2}{2} \right\} \right) \right|
		\left( \frac{1}{\varepsilon}\left(\frac{C}{\delta}\varepsilon + C\right)^p
		+ 3\varepsilon \left(\frac{C}{\delta^2} + \frac{C}{\varepsilon}\right)^2 \right) \\
	&\quad\quad+ 2\varepsilon \left(\mu_{\varepsilon,\delta} \left(
		\Omega_{0,2\delta} \backslash \Omega_{0,\delta} \times 
		\left(-\max \left\{ \frac{\varepsilon L^1}{2},
		\frac{\varepsilon L^2}{2} \right\}, \max \left\{ \frac{\varepsilon L^1}{2},
		\frac{\varepsilon L^2}{2} \right\}\right) \right)\right)^2 \\
	&\quad\leq {\mathcal H}^{N-1}(\Omega_{0,2\delta} \backslash \Omega_{0,\delta})
		\max \{ \varepsilon L^1, \varepsilon L^2 \}
		\left( \frac{1}{\varepsilon} \left(\frac{C}{\delta}\varepsilon + C \right)^p
		+ 3\varepsilon \left(\frac{C}{\delta^2} + \frac{C}{\varepsilon}\right)^2 \right)
		+ \varepsilon C.
	\end{align*}
	Choosing $\eta=\pm a\otimes e_N$ in Remark \ref{R6.1} we observe that
	\begin{equation*}
		W(\xi) \leq C |\min \{ \xi-a \otimes e_N, \xi + a \otimes e_N  \}|^p
	\end{equation*}
	which, together with the above estimates, implies that
	\begin{equation} \label{temp31}
	\limsup_{\delta \rightarrow 0^+} \limsup_{\varepsilon \rightarrow 0^+}
	E_\varepsilon(u_{\varepsilon,\delta}, \rho_{\varepsilon,\delta},
	\Omega \backslash \{ \psi_\delta = 1 \} \cap (\{ \varphi_\delta = 0 \}
	\cup \{ \varphi_\delta = 1 \})) = 0.
	\end{equation}
	We are left to estimate the energy on the set $\{ \psi_\delta = 1 \} \cap \left( \{ \varphi_\delta = 0 \}
	\cup \{ \varphi_\delta = 1 \} \right)$.
	To this end, without loss of generality we may suppose that
	$x_0 \notin \{ |x_N| \leq \frac{\varepsilon L^1}{2} \}$
	and $\tilde \delta \leq C$.
	We have that
	\begin{align*}
	E_\varepsilon(&u_{\varepsilon,\delta}, \rho_{\varepsilon,\delta},
		\{\psi_\delta = 1 \} \cap \{ \varphi_\delta = 1 \}) \leq
		\int_{K_\delta} \frac{1}{\varepsilon}
		W(\nabla u_{\varepsilon,\delta}) + \varepsilon 
		|\nabla^2 u_{\varepsilon, \delta}|^2 + \varepsilon
		(\rho_{\varepsilon,\delta} - |\nabla^2 u_{\varepsilon,\delta}|^2) \,dx \\
	&= \int_{K_\delta} \frac{1}{\varepsilon} W(\nabla z^1_\varepsilon)
		+ \varepsilon |\nabla^2 z^1_\varepsilon|^2 + \varepsilon
		(\rho_{\varepsilon,\delta} - |\nabla^2 z^1_\varepsilon|)^2 \,dx \\
	&= \int_{K_\delta \cap \{ |x_N| < \frac{\varepsilon L^1}{2} \}}
		\frac{1}{\varepsilon}W(\nabla z^1_\varepsilon)
		+ \varepsilon |\nabla^2 z^1_\varepsilon|^2 + \varepsilon
		(\max\{ |\nabla^2 z^1_\varepsilon| + \lambda^1_\varepsilon ,0 \} 
		- |\nabla^2 z^1_\varepsilon|)^2 \,dx \\
	&\quad+\int_{K_\delta \cap \{\frac{\varepsilon L^1}{2} < |x_N| 
		< \frac{\varepsilon L^1}{2} + \sqrt{\varepsilon} \}} \varepsilon
		\left( \frac{1}{2\sqrt{\varepsilon}} \left(\gamma_1 + \tilde \delta 
		- \int_{Q} \max\{|\nabla^2 v^1| + \lambda^1,0 \} \right) \right)^2 \,dx
		+ \sqrt{\varepsilon}\beta \chi_K(x_0) \\
	&\leq \int_{K_\delta \cap \{ |x_N| < \frac{\varepsilon L^1}{2} \}}
		\frac{1}{\varepsilon}W(\nabla z^1_\varepsilon)
		 + \varepsilon
		(\min\{ |\nabla^2 z^1_\varepsilon|^2 + (\lambda^1_\varepsilon)^2 ,
		2 |\nabla^2 z^1_\varepsilon|^2 \} \,dx  + C\sqrt{\varepsilon}\\
	&= \int_{K'_\delta \times(-\frac{1}{2},\frac{1}{2})} \hspace{-1cm}L^1
		W\left( \nabla v^1 \left( \frac{x'}{\varepsilon L^1},x_N \right) \right) + \frac{1}{L^1}
		\min \left\{ \left|\nabla^2 v^1 \left( \frac{x'}{\varepsilon L^1},x_N \right) \right|^2 + (\lambda^1)^2,
		2\left|\nabla^2 v^1 \left( \frac{x'}{\varepsilon L^1},x_N \right) \right|^2 \right\} \,dx \\
	&\quad+ C\sqrt{\varepsilon} .
		\end{align*}
Using the Riemann-Lebesgue lemma in the estimate above we get
\begin{equation*}
\lim_{\varepsilon \rightarrow 0^+}E_\varepsilon(u_{\varepsilon,\delta}, \rho_{\varepsilon,\delta},
		\{\psi_\delta = 1 \} \cap \{ \varphi_\delta = 1 \})\leq {\mathcal H}^{N-1}(K'_\delta)
		F_{1/L^1}(v^1,\lambda^1) 
		\leq {\mathcal H}^{N-1}(K'_\delta) (\Phi(\gamma_1) + \tilde \delta).
\end{equation*}

	An analogous argument leads to the estimate of the energy on the set
	$\{\psi_\delta = 1 \} \cap \{ \varphi_\delta = 0 \}$, namely
	\begin{equation*}
	\lim_{\varepsilon \rightarrow 0^+} E_\varepsilon
	(u_{\varepsilon,\delta}, \rho_{\varepsilon,\delta},
	\{\psi_\delta = 1 \} \cap \{ \varphi_\delta = 0 \})
	\leq {\mathcal H}^{N-1}(\Omega_{0,2\delta} \backslash K') (\Phi(\gamma_1) 
	+ \tilde \delta).
	\end{equation*}
	The last two inequalities together with (\ref{temp31}) imply that
	\begin{equation*}
	\lim_{\tilde \delta\rightarrow 0^+} 
	\lim_{\delta \rightarrow 0^+} \lim_{\varepsilon \rightarrow 0^+}
	E_\varepsilon (u_{\varepsilon,\delta}, \rho_{\varepsilon,\delta})\leq E(u,\mu).
	\end{equation*}
	Using a diagonal argument, we have shown the assertion of the theorem
	if $u$ and $\mu$ are as in \eqref{tempKorr7}.
\item[Step 1.2]
	We consider the case of finitely many interfaces, namely $S_{\nabla u}=\cup_{h=1}^n \Omega_{s_h}$, and
	\begin{equation*}
		\mu = \sum_{i=0}^n \gamma_i \chi_{K_i} + \sum_{i=0}^m \beta_i \delta_{x_i}, \quad
		\gamma_i, \beta_i \geq 0, K_i \subset S_{\nabla u} \text { compact and pairwise disjoint, } x_i 
	\in \Omega \backslash S_{\nabla u}.
	\end{equation*}
	By Theorem \ref{laminates}, we can suppose that near an interface $\Omega_{s_h}$, $u$ takes the form
	\begin{equation*}
		u(x) = \pm |x_n - s_h| + c_h.
	\end{equation*}
	Therefore, we can (up to adding constants) apply Step 1.1 to obtain a recovery sequence
	near the interface $\Omega_{s_h}$ for any summand in the definition of $\mu$.
	Since this construction is local and the sets $\Omega_{s_h}$ as well as 
	the sets $K_i$ and $\{ x_i \}$ have positive distance from each other,
	these local constructions can be glued to give a recovery sequence for $\mu$
	near $\Omega_{s_h}$ and then also for $u,\mu$ in $\Omega$.
	We leave the details to the reader.
\item[Step 1.3]
	We now consider the case of infinitely many interfaces, namely $S_{\nabla u}=\cup_{h=1}^\infty \Omega_{s_h}$ and we let $\mu$ be
	a finite sum of terms as in step 1.2.
	As in the proof  of \cite[Theorem 5.5 step 2]{cfl}, only $\alpha$ and $\beta$ can be accumulation points of the sequence $(s_h)$
	(otherwise one could find a cylinder $Z \subset \Omega$ with axis in direction of $e_N$
	that intersects infinitely
	many of the $\Omega_{s_h}$, and this would contradict the fact that 
	$\sum_h {\mathcal H}^{N-1}(\Omega_{s_h}) < \infty$).
	We can choose a decreasing sequence $\delta_k \rightarrow 0^+$
	such that $\{\alpha + \delta_k, \beta - \delta_k\}_k \cap \{s_h\}_h
	= \emptyset$ and such that $\text{supp} (\mu) \subset\subset 
	\{\alpha + \delta_k < x_N < \beta - \delta_k\}$ for all $k$.
	Then it follows that $u$ has only finitely many interfaces in
	\[U_k := \Omega \cap \{\alpha + \delta_k < x_N < \beta - \delta_k\}.\]
	Therefore, we can apply step 1.2 and find
	for any given sequence $\varepsilon_h \rightarrow 0^+$
	sequences $(u_h^k)_h$ and $(\mu_h)_h$ such that
	\begin{equation}\label{step1.3}
	\lim_{h \rightarrow \infty} u_h^k = u \text{ in } W^{1,1}(\Omega, U_k),
	\quad \mu_h \weak* \mu \text{ in } \Omega,
	\quad \lim_{h \rightarrow \infty} E_{\varepsilon_h}(u_h^k, \mu_h, U_k)
	\leq E(u,\mu,U_k),
	\end{equation}
	where we used that $\text{supp} (\mu) \subset\subset U_k$ and that
	therefore the sequences $\mu_h$ which approximate $\mu$ in Step 1.2
	can be chosen in such a way that they not depend on $k$.	
	By construction it holds
	\[E_{\varepsilon_h}(u_h^k,\mu_h,U_k) 
	= E_{\varepsilon_h}(u_h^k,\mu_h,\Omega)\]
	which implies by \eqref{step1.3} that
	\begin{equation} \label{temp32}
	\lim_{k \rightarrow \infty} \lim_{h \rightarrow \infty}
		E_{\varepsilon_h}(u_h^k,\mu_h,\Omega) 
		\leq \lim_{k \rightarrow \infty} E(u,\mu,U_k) \leq E(u,\mu,\Omega).
	\end{equation}
	The existence of  a subsequence $u_h^{k(h)}$ with the desired properties
	follows by a diagonal argument.
\item[Step 2]
	Now we assume that $\mu$ is a finite positive Radon measure.
	It holds
	\begin{equation*}
	0 \leq g := \frac{d \mu}{d \HN-1 \llcorner S_{\nabla u}} 
	\in L^1(S_{\nabla u},\HN-1).
	\end{equation*}
	We observe that there exist functions
	\begin{equation*}
	g_k = \sum_{i=0}^{n_k} \gamma_{k,i} \chi_{K'_{k,i}},
	\quad K'_{k,i} \subset S_{\nabla u} \text{ compact, pairwise disjoint}
	\end{equation*}
	which satisfy
	\begin{equation*}
	\lim_{k \rightarrow \infty} g_k = g \text{ in } L^1(S_{\nabla u},{\mathcal H}^{N-1}),
	\lim_{k \rightarrow \infty}{\mathcal H}^{N-1}(\{ g_k < g \}) = 0.
	\end{equation*}
	By the boundedness of $\mu$, we find measures $\sum_{i=0}^{m_k} \beta_{k,i} \delta_{x_{k,i}}$
	supported outside $S_{\nabla u}$ such that
	\begin{equation*}
	\sum_{i=0}^{m_k} \beta_{k,i} \delta_{x_{k,i}} \weak* \mu 
	- g {\mathcal H}^{N-1} \llcorner S_{\nabla u}.
	\end{equation*}
	It follows that
	\begin{equation*}
	\mu_k := g_k {\mathcal H}^{N-1} \llcorner S_{\nabla u}
	+ \sum_{i=0}^{m_k} \beta_{k,i} \delta_{x_{k,i}} \weak* \mu.
	\end{equation*}
	Since the $\gammalimsup$ is lower semicontinuous and $\Phi$
	is non-increasing and bounded, we can use the result of step 1 to obtain
	\begin{align*}
	\gammalimsup_{h \rightarrow \infty} E_{\varepsilon_h}(u,\mu)
		&\leq \liminf_{k \rightarrow \infty} 
		\gammalimsup_{h \rightarrow \infty} E_{\varepsilon_h}(u,\mu_k) \\
	&\leq \liminf_{k \rightarrow \infty} E(u,\mu_k)
		= \liminf_{k \rightarrow \infty} \int_\Omega 
		 \Phi(g_k) d {\mathcal H}^{N-1} \\
	&\leq \liminf_{k \rightarrow \infty} \int_\Omega \Phi(g) d {\mathcal H}^{N-1}
		+ C\lim_{k \rightarrow \infty}{\mathcal H}^{N-1}(\{ g_k < g \}) = E(u,\mu)
	\end{align*}
	for any sequence $\varepsilon_h \rightarrow 0^+$. This concludes the
	proof of the theorem.
\end{description}
\end{proof}

\end{document}